\documentclass[12pt]{amsart} 

\usepackage{amsmath} \usepackage{amssymb} \usepackage{amsthm}
\usepackage{amscd} 
\usepackage[all]{xy} 

\def\a{{\mathfrak{a}}} \def\b{{\mathfrak{b}}}\def\I{{\mathcal{I}}}
\def\F{{\mathbb{F}}} \def\J{{\mathcal{J}}} \def\m{{\mathfrak{m}}} \def\Z{{\mathbb{Z}}} \def\sO{{\mathcal{O}}}
\def\Q{{\mathbb{Q}}} \def\R{{\mathbb{R}}} 
\def\Hom{{\mathrm{Hom}}} \def\Ker{{\mathrm{Ker}}} 
  \def\Ann{{\mathrm{Ann}}}
\def\Spec{{\mathrm{Spec\; }}} 
\def\Div{{\mathrm{div}}} \def\adj{{\mathrm{adj}}}

\theoremstyle{plain}
\newtheorem{thm}{Theorem}[section] 
\newtheorem*{mainthm}{Main Theorem}
\newtheorem{cor}[thm]{Corollary}
\newtheorem{prop}[thm]{Proposition}

\newtheorem{propdef}[thm]{Proposition-Definition} 

\newtheorem{lem}[thm]{Lemma}
\theoremstyle{definition} 
\newtheorem{defn}[thm]{Definition}
\newtheorem{eg}[thm]{Example} 
\theoremstyle{remark}
\newtheorem{rem}[thm]{Remark}

\newtheorem*{claim}{Claim}
\newtheorem*{clproof}{Proof of Claim}
\newtheorem{cln}{Claim}
\newtheorem{clnproof}{Proof of Claim}
\newtheorem*{acknowledgement}{Acknowledgment}

\title{A characteristic $p$ analogue of\\ plt singularities and adjoint ideals}
\author{Shunsuke Takagi}
\address{Department of Mathematics, Kyushu University, 
6-10-1 Hakozaki, Higashi-ku, Fukuoka, 812-8581 JAPAN}
\email{stakagi@math.kyushu-u.ac.jp}

\begin{document}

\maketitle
\markboth{SHUNSUKE TAKAGI}{A characteristic $p$ analogue of adjoint ideals}

\begin{abstract}
We introduce a new variant of tight closure and give an interpretation of adjoint ideals via this tight closure. As a corollary, we prove that  a log pair $(X, \Delta)$ is plt if and only if the modulo $p$ reduction of $(X,\Delta)$ is divisorially F-regular for all large $p \gg 0$. 
Here, divisorially F-regular pairs are a class of singularities in positive characteristic introduced by Hara and Watanabe \cite{HW} in terms of Frobenius splitting. 
\end{abstract}

\section*{Introduction}
\label{intro}
The multiplier ideal sheaf $\J(X,D)$ associated to a log pair $(X,D)$ (i.e., $X$ is a normal complex variety and $D$ is an $\R$-divisor on $X$) is defined in terms of resolution of singularities and discrepancy divisors, and one can view this ideal sheaf as measuring how singular the pair $(X,D)$ is.  
However, when $X$ is smooth and $D$ is a (Cartier) integral divisor on $X$, the multiplier ideal sheaf $\J(X, D)$ is nothing but $\sO_X(-D)$ and does not reflect the singularities of $(X,D)$.
 On the other hand, the adjoint ideal sheaf $\adj(X,\Delta)$ of a boundary  $\Delta$ (i.e., $\Delta=\sum_i d_i\Delta_i$ is an $\R$-divisor with $0 \leq d_i \leq 1$) on $X$ is a variant of the multiplier ideal sheaf $\J(X,\Delta)$, and it encodes much information on the singularities of $(X,\Delta)$ even when $\Delta$ is a Cartier integral divisor. 
Ein-Lazarsfeld \cite{EL} and Debarre-Hacon \cite{DH} used the adjoint ideal sheaf to study the singularities of ample divisors of low degree on abelian varieties. 
Kawakita \cite{Ka} used the adjoint ideal sheaf to prove inversion of adjunction on log canonicity. 
The purpose of this paper is to give an interpretation of the adjoint ideal sheaf via a variant of tight closure. 

Tight closure is an operation defined on ideals or modules in positive characteristic. 
It was introduced by Hochster-Huneke \cite{HH1} in the 1980s. 
The notions of F-regular rings and F-rational rings are defined via tight closure, and they turned out to correspond to log terminal and rational singularities, respectively (\cite{Ha1}, \cite{HW}, \cite{MS}, \cite{Sm1}). 
This result is generalized to the correspondence of the test ideal and the multiplier ideal of the trivial divisor (\cite{Ha2}, \cite{Sm}). 
Here, the test ideal $\widetilde{\tau}(R)$ of a Noetherian local ring $(R,\m)$ of prime characteristic $p$ is the annihilator ideal of the tight closure $0^*_{E_R(R/\m)}$ of the zero submodule in the injective hull $E_R(R/\m)$ of the residue field $R/\m$ of $R$, and it plays a central role in the theory of  tight  closure.
Since we can enjoy the usefulness of multiplier ideals only when they are associated to various ideals or divisors, Hara-Yoshida \cite{HY} and the author \cite{Ta} introduced generalizations of tight closure and   of the test ideal associated to any given ideal and divisor, respectively.
They then proved the correspondence of their generalized test ideals and the multiplier ideals, building on earlier results of Hara \cite{Ha2} and Smith \cite{Sm}. 

In this paper, we introduce another generalization of tight closure associated to any given boundary, called divisorial tight closure, and investigate its properties. 
Let $(R,\m)$ be a Noetherian normal local ring of characteristic $p>0$ and $\Delta$ be a boundary on $X:=\Spec R$. Then the divisorial $\Delta$-tight closure $I^{\mathrm{div}*\Delta}$ of an ideal $I \subseteq R$ is the ideal consisting of all elements $x \in R$ for which there exists  $c \in R$ not in any minimal prime ideal of $H^0(X,\sO_X(-\lfloor \Delta \rfloor))$ such that 
$$cx^q \in I^{[q]}H^0(X,\sO_X((q-1)\Delta))$$
for all large $q=p^e$, where $I^{[q]}$ is the ideal generated by the $q$-th powers of elements of $I$. 
If $N \subseteq M$ are $R$-modules, then the divisorial $\Delta$-tight closure $N^{\mathrm{div}*\Delta}_M$ of $N$ in $M$ is defined similarly.  
We then define the divisorial test ideal $\widetilde{\tau}^{\mathrm{div}}(R,\Delta)$ to be the annihilator ideal of the divisorial $\Delta$-tight closure $0^{\mathrm{div}*\Delta}_{E_R(R/\m)}$ of the zero submodule in the injective hull $E_R(R/\m)$ of the residue field of $R$.  
In the case when $\Delta=0$, divisorial $\Delta$-tight closure coincides with classical tight closure and the divisorial test ideal $\widetilde{\tau}^{\mathrm{div}}(R,0)$ is nothing but the test ideal $\widetilde{\tau}(R)$. 
The divisorial test ideal $\widetilde{\tau}^{\mathrm{div}}(R,\Delta)$ establishes several nice properties analogous to those of the adjoint  ideal sheaf $\adj(X,\Delta)$: the restriction theorem (Theorem \ref{restriction}), the subadditivity theorem (Theorem \ref{subadditivity}), etc. 
By virtue of the following theorem, the divisorial test ideal $\widetilde{\tau}^{\mathrm{div}}(R,\Delta)$ can be viewed as a characteristic $p$ analogue of the adjoint ideal $\adj(X,\Delta)$. 
 
\begin{mainthm}[\textup{Theorem 5.3}]
Let $(R,\m)$ be a normal local ring essentially of finite type over a perfect field of characteristic $p>0$, and let $\Delta$ be a boundary on $X:=\Spec R$ such that $K_X+\Delta$ is $\R$-Cartier.    
Assume that $(R,\Delta)$ is reduced from characteristic zero to characteristic $p \gg 0$, together with a log resolution $f:\widetilde{X} \to X$ of $(X,\Delta)$ giving the adjoint ideal $\adj(X,\Delta)$. Then
$$\adj(X,\Delta)=\widetilde{\tau}^{\mathrm{div}}(R,\Delta).$$
\end{mainthm}

Plt singularities are one of the important classes of singularities arising in the minimal model program. 
Hara-Watanabe \cite{HW} introduced the notion of divisorially F-regular pairs in terms of Frobenius splitting, and they conjectured that divisorially F-regular pairs correspond to plt singularities.
The adjoint ideal sheaf $\adj(X,\Delta)$ defines the locus of non-plt points of $(X,\Delta)$ in $X$.  
Likewise, the divisorial test ideal $\widetilde{\tau}^{\mathrm{div}}(R,\Delta)$ defines the locus of non-divisorially-F-regular points of $(R,\Delta)$ in $\Spec R$.  
Thus, the conjecture of Hara and Watanabe immediately follows from Main Theorem. 

\begin{acknowledgement}
The author is indebted to Yasunari Nagai and Ryo Takahashi for valuable conversations and Ken-ichi Yoshida for useful comments. 
He also would like to thank Masayuki Kawakita and the referee for helpful suggestions. 
This work was partially supported by Grant-in-Aid for Young Scientists (B) of Japan Society for the Promotion of Science.
\end{acknowledgement}

\section{Multiplier ideals and Adjoint ideals}
\label{sec:1}
In this section, we briefly review the definition and local properties of multiplier ideals and adjoint ideals.  Our main reference is \cite{La}.  

Let $X$ be a normal algebraic variety over a field $k$ of characteristic zero.  For an $\R$-divisor $D=\sum_i r_i D_i$ on $X$, we define 
\begin{align*}
\lfloor D \rfloor&:=\sum_i \lfloor r_i \rfloor D_i  \; \; \textup{the \textit{round down} of $D$},\\
\lceil D \rceil&:=\sum_i \lceil r_i \rceil D_i=-\lfloor -D \rfloor \; \; \textup{the \textit{round up} of $D$},\\
\{ D \}&:=\sum_i \{ r_i \} D_i=D-\lfloor D \rfloor \; \; \textup{the \textit{fractional part} of $D$},
\end{align*}
where for $r \in \R$, $\lfloor r \rfloor$ is the largest integer less than or equal to $r$. 

Let $\Delta=\sum d_i \Delta_i$ be an effective $\R$-divisor on $X$ such that  $K_X+\Delta$ is $\R$-Cartier, and let $\a \subseteq \sO_X$ be an ideal sheaf. 
A \textit{log resolution} of $((X,\Delta);\a)$ is a proper birational morphism $f: \widetilde{X} \to X$ with $\widetilde{X}$ nonsingular such that $\a\sO_{\widetilde{X}}=\sO_{\widetilde{X}}(-F)$ is invertible and $\mathrm{Exc}(f) \cup \mathrm{Supp}(f^{-1}_*\Delta+F)$ is a simple normal crossing divisor. 
The existence of log resolutions is guaranteed by Hironaka's desingularization theorem \cite{Hi}. 
Let $f: \widetilde{X} \to X$ be a log resolution of $((X,\Delta);\a)$.
Then there are canonically defined real numbers $a(E)=a(E,X,\Delta)$, called the \textit{discrepancies} of $E$ with respect to the pair $(X,\Delta)$, attached to each prime divisor $E$ on $\widetilde{X}$ having the property that
$$K_{\widetilde{X}} \equiv_{\mathrm{num}} f^*(K_X+\Delta)+\sum a(E,X,\Delta) E,$$
the sum running over all prime divisors $E$ on $\widetilde{X}$
(we define $a(f^{-1}_*\Delta_i, X, \Delta)=-d_i$ and $a(E,X,\Delta)=0$ for prime divisors $E$ not supported on $\mathrm{Exc}(f) \cup \mathrm{Supp}(f^{-1}_*\Delta)$).  
Note that $a(E,X,\Delta)$ is independent of $f$ when one views $E$ as being a valuation of the function field $k(X)$ of $X$. 

\begin{defn}[\textup{\cite[Definition 9.3.60]{La}}]
In the above situation, let $t>0$ be a real number. 
Fix a log resolution $f: \widetilde{X} \to X$ of $((X,\Delta);\a)$ so that $\a\sO_{\widetilde{X}}=\sO_{\widetilde{X}}(-F)$ for an effective divisor $F=\sum b(E) E$ on $\widetilde{X}$. 
The \textit{multiplier ideal sheaf} $\J((X,\Delta);\a^t)$ of $\a$ with exponent $t$ for the pair $(X,\Delta)$ is
$$\J((X,\Delta);\a^t)=f_*\sO_{\widetilde{X}}(K_{\widetilde{X}}-\lfloor f^*(K_X+\Delta)+tF \rfloor)=f_*\sO_{\widetilde{X}}(\sum \lceil a(E)-t \cdot b(E) \rceil E).$$
We denote this ideal simply by $\J(X,\Delta)$ (resp. $\J(X,\a^t)$) when $\a=\sO_X$ (resp. $\Delta=0$).   
$\J((X,\Delta);\a^t)$ is independent of the choice of the log resolution $f$ (\cite[Theorem 9.2.18]{La}). 
\end{defn}

Next we define a variant of the multiplier ideal sheaf, called the adjoint ideal sheaf\footnote{
 Lipman called the multiplier ideal associated to an ideal $\a$ in a regular ring as the adjoint ideal of $\a$: however we use this term for a somewhat different notion according to \cite{La}.}. 
Kawakita \cite{Ka} used this ideal sheaf to prove inversion of adjunction on log canonicity and Shokurov \cite{Sh} also used it implicitly. 
 
An $\R$-divisor $\Delta=\sum_i d_i \Delta_i$ on $X$ is called a \textit{boundary} if $0 \leq d_i \leq 1$ for every $i$. 
Let $\Delta$ be a boundary on $X$ such that  $K_X+\Delta$ is $\R$-Cartier, and let $\a \subseteq \sO_X$ be an ideal sheaf such that no component of $\lfloor \Delta \rfloor$ is contained in the zero-locus of $\a$.

\begin{defn}[\textup{cf.\cite[Definition 9.3.47]{La}}]\label{adjoint}
In the above situation, let $t>0$ be a real number. 
Fix a log resolution $f: \widetilde{X} \to X$ of $((X,\Delta);\a)$ such that $\a\sO_{\widetilde{X}}=\sO_{\widetilde{X}}(-F)$ for an effective divisor $F$ on $\widetilde{X}$ and $f^{-1}_*{\lfloor \Delta \rfloor}$ is nonsingular (but possibly disconnected). 
Then the \textit{adjoint ideal sheaf} $\adj((X,\Delta);\a^t)$ of $\a$ with exponent $t$ for the pair $(X,\Delta)$ is
$$\adj((X,\Delta);\a^t)=f_*\sO_{\widetilde{X}}(K_{\widetilde{X}}-\lfloor f^*(K_X+\Delta)+tF \rfloor+f^{-1}_*\lfloor \Delta \rfloor).$$
We denote this ideal simply by $\adj(X,\Delta)$ when $\a=\sO_X$.  
\end{defn}

The multiplier ideal sheaf $\J(X,D)$ associated to an integral divisor $D$ on a smooth variety $X$ is just $\sO_X(-D)$ (\cite[Example 9.2.12]{La}). On the other hand, the adjoint ideal sheaf $\adj(X,D)$ associated to $D$ does carry significant information about the singularities of $D$.

There are important classes of singularities defined in terms of the triviality of the multiplier ideal sheaf or the adjoint ideal sheaf. 
\begin{defn}\label{pair}
Let $X$ be a normal variety over a field $k$ of characteristic zero, and let $\Delta$ be an effective $\R$-divisor on $X$ such that  $K_X+\Delta$ is $\R$-Cartier. 
\renewcommand{\labelenumi}{(\roman{enumi})}
\begin{enumerate}
\item We say that the pair $(X,\Delta)$ is Kawamata log terminal (or klt for short) if $\J(X,\Delta)=\sO_X$. 
\item We say that the pair $(X,\Delta)$ is purely log terminal (or plt for short) if $\Delta$ is a boundary (i.e., $\Delta=\sum_i d_i \Delta_i$ with $0 \leq d_i \leq 1$) and  $\adj(X,\Delta)=\sO_X$.  
\end{enumerate}
\end{defn}

\begin{rem}\label{adjoint rem}
\renewcommand{\labelenumi}{(\roman{enumi})}
\begin{enumerate}
\item One can verify by an argument similar to \cite[Theorem 9.2.18]{La} that the adjoint ideal sheaf $\adj((X,\Delta);\a^t)$ is independent of the choice of the log resolution used to define it. 

\item Obviously $\J((X,\Delta);\a^t) \subseteq \adj((X,\Delta);\a^t)$ (when $\Delta$ is a boundary). If $\lfloor \Delta \rfloor=0$, then the both ideal sheaves coincide with each other. 

\item Klt singularities are plt. 

\item We can define the adjoint ideal sheaf in more general setting. 
Let $D$ be a reduced integral divisor and $B$ be an effective $\R$-divisor on $X$ such that $D$ and $B$ have no common components and that $K_X+D+B$ is $\R$-Cartier.  
Let $\a \subseteq \sO_X$ be an ideal sheaf such that no component of $D$ is  contained in the zero-locus of $\a$. Then for any real number $t>0$, the adjoint ideal sheaf $\adj((X,D;B);\a^t)$ is defined to be
$$\adj((X,D;B);\a^t)=f_*\sO_{\widetilde{X}}(K_{\widetilde{X}}-\lfloor f^*(K_X+D+B)+tF \rfloor+f^{-1}_*D).$$
To avoid being too technical, we do not deal with this ideal in this paper. 
However, all results on the adjoint ideal sheaf $\adj((X,\Delta);\a^t)$ are extendable to those on this adjoint ideal sheaf $\adj((X,D;B);\a^t)$. See also Remark \ref{div t.c. remark} (iii). 
\end{enumerate}
\end{rem}

One of the most important local properties of adjoint ideal sheaves is the restriction theorem. 
\begin{thm}[\textup{cf. \cite[Theorem 9.5.16]{La}}]\label{restriction in char 0}
Let $X$ be a normal variety over a field of  characteristic zero and $\Delta$ be a boundary on $X$ such that $K_X+\Delta$ is $\R$-Cartier. Put $D=\lfloor \Delta \rfloor$ and $B=\Delta-\lfloor \Delta \rfloor$.  
Let $\nu:D^{\nu} \to D$ be the normalization of $D$ and $B^{\nu}$ be the different of $B$ on $D^{\nu}$. That is, $B^{\nu}$ is an effective $\R$-divisor on $D^{\nu}$ such that $K_{D^{\nu}}+B^{\nu}=\nu^*((K_X+\Delta)|_D)$ $($the reader is referred to \cite[\S 3]{Sh} for details$)$. 
Let $\a \subseteq \sO_X$ be an ideal sheaf such that no component of $D$ is  contained in the zero-locus of $\a$. 
Then for any real number $t>0$, the sheaf $\nu_*\J((D^{\nu},B^{\nu});{\a\sO_{D^{\nu}}}^t)$ is an ideal sheaf of $D$ and one has 
$$\nu_*\J((D^{\nu},B^{\nu});{\a \sO_{D^{\nu}}}^t)= \adj((X,\Delta); \a^t) \sO_D.$$
In particular, if the adjoint ideal sheaf $\adj(X,\Delta)$ is trivial, then $D=\lfloor \Delta \rfloor$ is normal and Cohen-Macaulay. 
\end{thm}

The above restriction theorem is used to prove Demailly-Ein-Lazarsfeld's subadditivity property of multiplier ideals. 
\begin{thm}[\textup{cf. \cite{DEL}, \cite[Theorem 9.5.20]{La}}]\label{sub in char 0}
Let $X$ be a $d$-dimensional nonsingular variety over a field of characteristic zero. 
Fix a point $x \in X$ and  denote by $\m=\m_x$ the maximal ideal sheaf at $x$. 
If $\a,\b \subseteq \sO_X$ are ideal sheaves, then
$$\J(X,\a^s\b^t)_x \subseteq \sum_{\lambda+\mu=d} \J(X,\m^{\lambda}\a^s)_x\J(X,\m^{\mu}\b^t)_x$$
for any real numbers $s, t>0$. 
\end{thm}

\section{F-singularities of pairs and generalized test ideals}
\label{sec:2}
In this section, we recall the definition of F-singularities of pairs and $(\Delta, \a^t)$-tight closure used to define generalized test ideals. 
The reader is referred to \cite{HW}, \cite{HY} and \cite{Ta} for details. 

Throughout this paper, all rings are excellent reduced Noetherian commutative rings with identity. 
Let $R$ be a normal domain with quotient field $K$ and $D$ be an $\R$-divisor (not necessarily effective) on $\Spec R$. 
We denote 
$$R(D):=H^0(X,\sO_X(D))=\{0\} \cup \{x \in K \mid \Div_X(x)+D \geq 0 \}.$$
Let $\Delta$ be an effective $\R$-divisor on $\Spec R$, and denote by $R^{\circ, \Delta}$ the set of elements of $R$ which are not in any minimal prime ideal of $R(-\lfloor \Delta \rfloor) \subseteq R$. 
We write this set simply by $R^{\circ}$ when $\Delta=0$. 

In addition, suppose that $R$ is of characteristic $p>0$. 
For each $q=p^e$, $R((q-1)\Delta)$ viewed as an $R$-module via the $e$-times Frobenius map $F^e \colon R \to R((q-1)\Delta)$ sending $x$ to $x^q$ is denoted by ${}^e\! R((q-1)\Delta)$.
Since $R((q-1)\Delta)$ is a submodule of $K$, we can identify $F^e:R \to {}^e\! R((q-1)\Delta)$ with the natural inclusion map $R \hookrightarrow R((q-1)\Delta)^{1/q}$.
We say that $R$ is {\it F-finite} if ${}^1\! R$ (or $R^{1/p}$) is a finitely generated $R$-module. 
For example, any algebra essentially of finite type over a perfect field is F-finite. 
Also, for any ideal $I$ of $R$, we denote by $I^{[q]}$ the ideal of $R$ generated by the $q$-th powers of elements of $I$. 

\begin{defn}[\textup{[HW, Definition 2.1]}]\label{F-pair}
Let $R$ be an F-finite normal domain of characteristic $p>0$ and $\Delta$ be an effective $\R$-divisor on $\Spec R$.
\renewcommand{\labelenumi}{(\roman{enumi})}
\begin{enumerate}
\item $(R , \Delta)$ is said to be \textit{strongly F-regular} if for every $c \in R^{\circ}$, there exists $q=p^e$ such that $c^{1/q}R \hookrightarrow R((q-1)\Delta)^{1/q} $ splits as an $R$-module homomorphism.

\item $(R,\Delta)$ is said to be \textit{divisorially F-regular}\footnote{As we state in Corollary \ref{plt}, divisorially F-regular pairs correspond to plt singularities, not to dlt singularities. The reader is referred to \cite{KM} for the difference between plt singularities and dlt singularities.  Our term might be misleading in this sense, but  we use the same term as that in the original paper \cite{HW}.} if for every $c \in R^{\circ, \Delta}$, there exists $q=p^e$ such that $c^{1/q}R \hookrightarrow R((q-1)\Delta)^{1/q}$ splits as an $R$-module homomorphism. 
\end{enumerate}
\end{defn}

Let $R$ be a normal domain of characteristic $p > 0$ and $\Delta$ be an effective $\R$-divisor on $\Spec R$. Let $M$ be an $R$-module. 
For each $q=p^e$, we denote $\F^{e,\Delta}(M)= \F_{R}^{e,\Delta}(M) := {}^e\! R((q-1)\Delta) \otimes_R M$ and 
regard it as an $R$-module by the action of $R$ from the left. 
Then we have the $e$-times  Frobenius map $F_M^{e} \colon M \to 
\F^{e,\Delta}(M)$ induced on $M$. The image of an element $z \in M$ via this map is denoted by $z^q:= F_M^{e}(z) 
\in \F^{e,\Delta}(M)$. 
For an $R$-submodule $N$ of $M$, we denote by $N^{[q],\Delta}_{M}$ the 
image of the induced map $\F^{e,\Delta}(N) \to \F^{e,\Delta}(M)$. 
If $I $ is an ideal of $R$, then $I^{[q],\Delta}_{R}=I^{[q]}R((q-1)\Delta)$. 

\begin{defn}[\textup{cf. \cite[Definition 2.1]{Ta}, \cite[Definition 6.1]{HY}}]\label{tight closure}
In the above situation, let $\a$ be a nonzero ideal of $R$ and $t>0$ be a real number. 
\renewcommand{\labelenumi}{(\roman{enumi})}
\begin{enumerate}
\item If $N \subseteq M$ are $R$-modules, then the \textit{$(\Delta,\a^t)$-tight closure} $N^{*(\Delta,\a^t)}_M$ of $N$ in $M$ is defined to be the submodule of $M$ consisting of all elements $z \in M$ for which there exists $c \in R^{\circ}$ such that 
$$c\a^{\lceil tq \rceil}z^q \subseteq N^{[q],\Delta}_M$$ 
for all large $q = p^e$.

\item Let $E=\oplus_{\m} E_R(R/\m)$ be the direct sum, taken over all maximal ideals $\m$ of $R$, of the injective hulls of the residue fields $R/\m$. 
Then we define the \textit{generalized test ideal} $\widetilde{\tau}((R,\Delta);\a^t)$ by 
$$\widetilde{\tau}((R,\Delta);\a^t)=\mathrm{Ann}_R(0^{*(\Delta,\a^t)}_E) \subseteq R. $$ 
We denote this ideal simply by $\widetilde{\tau}(R,\Delta)$ (resp. $\widetilde{\tau}(R,\a^t)$) when $\a=R$ (resp. $\Delta=0$). 
\end{enumerate}
\end{defn}

\begin{rem}\label{F-pair remark}
\renewcommand{\labelenumi}{(\roman{enumi})}
\begin{enumerate}
\item $R$ is a strongly F-regular ring if and only if the pair $(R,0)$ is divisorially F-regular, or equivalently strongly F-regular. We refer the reader to \cite{HH1}, \cite{HH2} for strongly F-regular rings. 

\item Strongly F-regular pairs are divisorially F-regular.  When $\lfloor \Delta \rfloor=0$, the converse implication also holds true. 

\item $(\textup{\cite[Proposition 2.2 (3), (4)]{HW}})$ If the pair $(R,\Delta)$ is strongly F-regular (resp. divisorially F-regular), then $\lfloor \Delta \rfloor=0$ (resp. $\Delta$ is a boundary). 

\item (cf. $\!\textup{\cite[Proposition 2.2 (2)]{HW}}$)  
The definition of  divisorial F-regularity (resp. strong F-regularity) does not change even if we replace $R((q-1)\Delta)^{1/q}$ by $R(\lceil q \Delta \rceil-\lfloor \Delta \rfloor)^{1/q}$ (resp. $R(\lceil q \Delta \rceil)^{1/q}$) in Definition \ref{F-pair}. 
More generally, we can prove that $(R,\Delta)$ is divisorially F-regular (resp. strongly F-regular) if and only if for every $c \in R^{\circ, \Delta}$ (resp. $c \in R^{\circ}$), there exists $q'$ such that $c^{1/q}R \hookrightarrow R(\lceil q \Delta \rceil-\lfloor \Delta \rfloor)^{1/q}$ (resp. $c^{1/q}R \hookrightarrow R(\lceil q \Delta \rceil)^{1/q}$) splits as an $R$-module homomorphism for all  $q=p^e \geq q'$. 

\item $(\textup{\cite[Lemma 2.3]{Ta}})$ $(R,\Delta)$ is strongly F-regular if and only if $\widetilde{\tau}(R,\Delta)=R$. 
\end{enumerate}
\end{rem}

In order to relate the generalized ideal $\widetilde{\tau}((R,\Delta);\a^t)$ with the multiplier ideal $\J((X,\Delta);\a^t)$, we use the technique of ``reduction to characteristic $p \gg 0$." 

Let $R$ be an algebra essentially of finite type over a field $k$ of characteristic zero and $\Delta$ be an effective $\R$-divisor on $\Spec R$. 
Let $\a \subseteq  R$ be an ideal such that $\a \cap R^{\circ} \ne \emptyset$ and $t>0$ be a real number. 
One can choose a finitely generated $\Z$-subalgebra $A$ of $k$ and a subalgebra $R_A$ of $R$ essentially of finite type over $A$ such that the natural map $R_A \otimes_A k \to R$ is an isomorphism, 
$\rho^*\Delta_A=\Delta$ and ${\a_A}R=\a$, where $\rho:\Spec R \to \Spec R_A$ is the map associated to the inclusion $R_A \hookrightarrow R$, $\Delta_A:=\rho_*\Delta$ and ${\a_A}:=\a \cap R_A \subseteq R_A$. 
Given a closed point $s \in \Spec A$ with residue field $\kappa=\kappa(s)$, we denote the corresponding fibers over $s$ by $R_{\kappa}, \Delta_{\kappa}, {\a_{\kappa}}$. 
Then we refer to such $(\kappa, R_{\kappa}, \Delta_{\kappa}, {\a_{\kappa}})$ for a general closed point $s \in \Spec A$ with residue field $\kappa=\kappa(s)$ of sufficiently large characteristic $p \gg 0$ as ``{\it reduction to characteristic $p \gg 0$}'' of $(k,R, \Delta, \a)$, and the triple $(R_{\kappa},  \Delta_{\kappa}, {\a_{\kappa}}^{t})$ inherits the properties possessed by the original one $(R, \Delta, \a^{t})$ (how large $p$ has to be is depending on $t$).
Furthermore, given a log resolution $f:\widetilde{X} \to X=\Spec R$ of $((X, \Delta); \a)$, we can reduce this entire setup to characteristic $p \gg 0$. 

\begin{defn}
In the above situation, 
$(R,\Delta)$ is said to be of \textit{strongly F-regular} (resp. \textit{divisorially F-regular}) \textit{type} if reduction to characteristic $p \gg 0$ of $(R,\Delta)$ is strongly F-regular (resp. divisorially F-regular). 
\end{defn}

The multiplier ideal $\J((X,\Delta);\a^t)$ coincides, after reduction to characteristic $p \gg 0$, with the generalized ideal $\widetilde{\tau}((R,\Delta);\a^t)$. 
\begin{thm}[\textup{cf. \cite[Theorem 3.2]{Ta}, \cite[Theorem 6.8]{HY}}]\label{tau=mult}
Let $(R,\m)$ be a normal local ring essentially of finite type over a perfect field of positive characteristic $p$, and let $\Delta$ be an effective $\R$-divisor on $X:=\Spec R$ such that $K_X+\Delta$ is $\R$-Cartier. 
Let $\a \subseteq R$ be a nonzero ideal and $t >0$ be a fixed real number.   
Assume that $((R,\Delta);\a)$ is reduced from characteristic zero to characteristic $p \gg 0$, together with a log resolution $f:\widetilde{X} \to X$ of $((X,\Delta);\a)$ giving the multiplier ideal $\J((X,\Delta);\a^t)$. Then
$$\J((X,\Delta);\a^t)=\widetilde{\tau}((R,\Delta);\a^t).$$
In particular, $(X,\Delta)$ is klt if and only if $(R,\Delta)$ is of strongly F-regular type. 
\end{thm}

\section{Divisorial tight closure}
\label{sec:3}
In this section, we introduce another variant of tight closure, called divisorial tight closure, and investigate its basic properties. 
\begin{defn}\label{div tight closure}
Let $R$ be an F-finite normal domain of characteristic $p>0$ and $\Delta$ be a boundary (i.e., $\Delta=\sum_i d_i \Delta_i$ is an $\R$-divisor with $0 \leq d_i \leq 1$) on $\Spec R$. 
Let $\a \subseteq R$ be an ideal such that $\a \cap R^{\circ,\Delta} \ne \emptyset$ where 
$R^{\circ, \Delta}$ is the set of elements of $R$ which are not in any minimal prime ideal of $R(-\lfloor \Delta \rfloor) \subseteq R$, and let $t>0$ be a fixed real number. 
\renewcommand{\labelenumi}{(\roman{enumi})}
\begin{enumerate}
\item If $N \subseteq M$ are (not necessarily finitely generated) $R$-modules, then the \textit{divisorial $(\Delta,\a^t)$-tight closure} $N^{\mathrm{div}*(\Delta,\a^t)}_M$ of $N$ in $M$ is defined to be the submodule of $M$ consisting of all elements $z \in M$ for which there exists $c \in R^{\circ,\Delta}$ such that 
$$c\a^{\lceil tq \rceil}z^q \subseteq N^{[q],\Delta}_{M}$$ 
for all large $q = p^e$ (see the paragraph above Definition \ref{tight closure} for the meaning of the notation $N^{[q],\Delta}_{M}$). 
The divisorial $(\Delta,\a^t)$-tight closure of an ideal $I \subseteq R$ is defined by $I^{\mathrm{div}*(\Delta,\a^t)} := I^{\mathrm{div}*(\Delta,\a^t)}_R$. 
That is, $x \in R$ is in $I^{\mathrm{div}*(\Delta,\a^t)}$ if and only if there exists $c \in R^{\circ,\Delta}$ such that $cx^q \a^{\lceil tq \rceil} \subseteq I^{[q]}R((q-1)\Delta)$ for all large $q=p^e$. 
We denote $N^{\mathrm{div}*(\Delta,\a^t)}_M$ (resp. $I^{\mathrm{div}*(\Delta, \a^t)}$) simply by $N_M^{\mathrm{div}*\Delta}$ (resp. $I^{\mathrm{div}*\Delta}$) when $\a=R$. 

\item Let $E=\oplus_{\m} E_R(R/\m)$ be the direct sum, taken over all maximal ideals $\m$ of $R$, of the injective hulls of the residue fields $R/\m$. 
We define the \textit{divisorial test ideal} $\widetilde{\tau}^{\mathrm{div}}((R,\Delta);\a^t)$ by 
$$\widetilde{\tau}^{\mathrm{div}}((R,\Delta);\a^t)=\mathrm{Ann}_R(0^{\mathrm{div}*(\Delta,\a^t)}_E) \subseteq R. $$
We denote this ideal simply by $\widetilde{\tau}^{\mathrm{div}}(R,\Delta)$ when $\a=R$. 
\end{enumerate}
\end{defn}

\begin{eg}
\renewcommand{\labelenumi}{(\roman{enumi})}
\begin{enumerate}
\item When $\Delta=\Div_R(f)$ is an irreducible Cartier divisor, an element $x \in R$ is in $I^{\mathrm{div}*(\Delta,\a^t)}$ if and only if there exists $c \in R \setminus fR$ such that $cf^{q-1}\a^{\lceil tq \rceil}x^q \in I^{[q]}$ for all large $q=p^e$. 
\item Let $R=\F_p[[x,y]]$ be a complete regular local ring of characteristic $p>0$, $\Delta=\Div_R(xy)$ and $I=(x^2,y^2) \subset R$. 
Then 
$$I^{\mathrm{div}*\Delta}=(x^2,xy,y^2) \subsetneq (x,y)=(I^{\mathrm{div}*\Delta})^{\mathrm{div}*\Delta}.$$
Thus,  divisorial tight closure is not a closure operation in general. 
\end{enumerate}
\end{eg}

\begin{rem}\label{div t.c. remark}
\renewcommand{\labelenumi}{(\roman{enumi})}
\begin{enumerate}
\item The definition of divisorial tight closure extends to several exponents. 
Given ideals $\a_1, \dots, \a_r$ of $R$ with $\a_i \cap R^{\circ, \Delta} \ne \emptyset$ and real numbers $t_1, \dots, t_r >0$, 
one can define divisorial $(\Delta,\a_1^{t_1} \dots \a_r^{t_r})$-tight closure as follows: 
if $N \subseteq M$ are $R$-modules, then $z \in M$ is in the divisorial $(\Delta,\a_1^{t_1} \dots \a_r^{t_r})$-tight closure $N_M^{\mathrm{div}*(\Delta,\a_1^{t_1} \dots \a_r^{t_r})}$ of $N$ in $M$ if and only if there exists $c \in R^{\circ,\Delta}$ such that $c\a_1^{\lceil t_1q \rceil} \dots \a_r^{\lceil t_rq \rceil}z^q \subseteq N^{[q],\Delta}_{M}$ for all large $q=p^e$. 

\item As a generalization of the test ideal $\tau(R)$ in the classical tight closure theory, we can define another ideal ${\tau}^{\mathrm{div}}((R,\Delta);\a^t)$ by 
$${\tau}^{\mathrm{div}}((R,\Delta);\a^t)=\bigcap_{I \subseteq R}(I:I^{\mathrm{div}*(\Delta,\a^t)}), $$
where $I$ runs through all ideals $I$ of $R$. 
Then one can expect that this ideal coincides with the divisorial test ideal $\widetilde{\tau}^{\mathrm{div}}((R,\Delta);\a^t)$ under the assumption that $K_X+\Delta$ is $\R$-Cartier, as is a generalization of \cite[Theorem 1.13]{HY}, \cite[Lemma 3.4]{Sm} and \cite[Theorem 2.8]{Ta}; see also \cite{AM}.
We can prove it if there exists an integer $r$ not divisible by $p$ such that $r(K_X+\Delta)$ is Cartier, but it is a open problem in general.  
In addition, compared with $\widetilde{\tau}^{\mathrm{div}}((R,\Delta);\a^t)$, the ideal ${\tau}^{\mathrm{div}}((R,\Delta);\a^t)$ seems to be more difficult to handle. 
These are the reason why we don't pursue properties of the ideal ${\tau}^{\mathrm{div}}((R,\Delta);\a^t)$ in this paper. 

\item We can define the divisorial test ideal in more general setting. 
Let $D$ be a reduced integral divisor and $B$ be an effective $\R$-divisor on $\Spec R$ such that $D$ and $B$ have no common components. Let $\a \subseteq R$ be an ideal such that $\a \cap R^{\circ, D} \ne \emptyset$ and $t>0$ be a real number. 
The divisorial $(D;B,\a^t)$-tight closure $0^{\mathrm{div}*(D;B,\a^t)}_E$ of the zero submodule in $E$ is defined to be the submodule of $E$ consisting of all elements $z \in E$ for which there exists $c \in R^{\circ,D}$ such that 
$$c\a^{\lceil tq \rceil}z^q \subseteq 0^{[q],D+B}_{E}$$ 
for all large $q = p^e$. 
We then define the divisorial test ideal $\widetilde{\tau}^{\mathrm{div}}((R,D;B);\a^t)$ by 
$$\widetilde{\tau}^{\mathrm{div}}((R,D;B);\a^t)=\mathrm{Ann}_R(0^{\mathrm{div}*(D;B,\a^t)}_E). $$
This divisorial test ideal $\widetilde{\tau}^{\mathrm{div}}((R,D;B);\a^t)$ corresponds to the adjoint ideal 
$\adj((X,D;B);\a^t)$ defined in Remark \ref{adjoint rem} (iv). 
\end{enumerate}
\end{rem}

We list up basic properties of the divisorial test ideal $\widetilde{\tau}^{\mathrm{div}}((R,\Delta);\a^t)$ in the following. 
\begin{prop}\label{basic}
Let the notation be as in Definition $\ref{div tight closure}$.
\renewcommand{\labelenumi}{$(\arabic{enumi})$}
\begin{enumerate}
\item For any effective $\R$-divisor $\Delta' \leq \Delta$ on $\Spec R$, for any ideal $\a \subseteq \b \subseteq R$ and for any real number $s \leq t$,  one has 
$$\widetilde{\tau}^{\mathrm{div}}((R,\Delta);\a^t) \subseteq \widetilde{\tau}^{\mathrm{div}}((R,\Delta');\b^s).$$ 

\item Let $\Delta'$ be another boundary on $\Spec R$ such that $\Delta+\Delta'$ is also a boundary and $\a \cap R^{\circ, \Delta'} \ne \emptyset$. 
Let $\b \subseteq R$ be an ideal such that  $\b \cap R^{\circ, \Delta+\Delta'} \ne \emptyset$. 
Then 
$$R(-\Delta') \widetilde{\tau}^{\mathrm{div}}((R,\Delta);\a^t) \b \subseteq \widetilde{\tau}^{\mathrm{div}}((R,\Delta+\Delta');\a^t\b).$$

\item 
$(R,\Delta)$ is divisorially F-regular if and only if $\widetilde{\tau}^{\mathrm{div}}(R, \Delta)=R$. 
\end{enumerate}
\end{prop}

\begin{proof}
$(1)$ The proof is similar to that of \cite[Proposition 2.2 (5)]{HW}. 
The non-trivial case is only when $\mathrm{Supp}(\lfloor \Delta \rfloor)\setminus \mathrm{Supp}(\lfloor \Delta' \rfloor) \ne \emptyset$. 
To prove the assertion in this case, we may assume that without loss of generality that there exists a unique irreducible component $\Delta_0$ of $\lfloor \Delta \rfloor$ such that $\Delta_0 \not\subseteq \mathrm{Supp}(\lfloor \Delta' \rfloor)$. 
Let $z \in 0_E^{\mathrm{div}*(\Delta',\b^s)}$. By definition, there exists $c \in R^{\circ, \Delta'}$ such that 
$c\b^{\lceil sq \rceil}z^q=0$ in $\F^{e,\Delta'}(E)$ for all large $q=p^e$. 
We may assume that $c$ is in $\mathfrak{p}=R(-\Delta_0)$, and let $\nu=v_{\mathfrak{p}}(c)$ be the value of $c$ at $\mathfrak{p}$. Then one can choose $d \in R^{\circ, \Delta}$ which is in $cR(\nu \Delta_0)$. 
This implies that $d R((q-1)\Delta') \subseteq c R((q-1)\Delta)$ for all large $q=p^e$. 
Since $\a^{\lceil tq \rceil} \subseteq \b^{\lceil sq \rceil}$, we have $d\a^{\lceil tq \rceil}z^q=0$ in $\F^{e,\Delta}(E)$ for all large $q=p^e$, that is, $z \in 0_E^{\mathrm{div}*(\Delta,\a^t)}$.

$(2)$ It is enough to show that $0_E^{\mathrm{div}*(\Delta+\Delta',\a^t\b)} \subseteq (0_E^{\mathrm{div}*(\Delta,\a^t)}:R(-\Delta')\b)_E$. 
Let $z \in 0_E^{\mathrm{div}*(\Delta+\Delta',\a^t\b)}$. 
By definition, there is $c \in R^{\circ,\Delta+\Delta'}$ such that $c\a^{\lceil tq \rceil}\b^qz^q=0$ in $\F^{e, \Delta+\Delta'}(E)$ for all large $q=p^e$. 
Since $R(-\Delta')^{[q]} R((q-1)(\Delta+\Delta')) \subseteq R((q-1)\Delta)$ and $\b^{[q]} \subseteq \b^q$, one has $c\a^{\lceil tq \rceil}R(-\Delta')^{[q]}\b^{[q]}z^q=0$ in $\F^{e,\Delta}(E)$. 
This implies that $z$ lies in $(0_E^{\mathrm{div}*(\Delta,\a^t)}:R(-\Delta')\b)_E$. 

$(3)$ The assertion follows from Remark \ref{F-pair remark} and Lemma \ref{key lemma}. 
\end{proof}

\begin{lem}\label{key lemma}
Let $(R,\m)$ be an F-finite normal local ring of characteristic $p > 0$ and $\Delta$ be a boundary on $\Spec R$. 
Let $\a \subseteq R$ be an ideal such that $\a \cap R^{\circ,\Delta} \ne \emptyset$ and $t > 0$ be a real number. 
Fix a system of generators $x^{(e)}_1,\ldots, x^{(e)}_{r_{e}}$ of $\a^{\lceil tq \rceil}$ for each $q = p^e$. Then an element $c \in R$ lies in $\widetilde{\tau}^{\mathrm{div}}((R,\Delta);\a^t)$ if and only if for any $d \in R^{\circ,\Delta}$ and any positive integer $e_0$, there exist an integer $e_1 \geq e_0$ and $R$-homomorphisms $\varphi^{(e)}_i \in \Hom_R(R(\lceil p^e \Delta \rceil-\lfloor \Delta \rfloor)^{1/p^e},R)$ for $e_0 \leq e \leq e_1$ and $1 \leq i \leq r_e$ such that 
$$c=\sum_{e=e_0}^{e_1}\sum_{i=1}^{r_e} \varphi^{(e)}_i((dx_i^{(e)})^{1/p^e}).$$
\end{lem}
\begin{proof}
The proof is identical to that for the generalized test ideal $\widetilde{\tau}(R,\a^t)$. See \cite[Lemma 2.1]{HT}.  
\end{proof}

As immediate applications of Lemma \ref{key lemma}, we can show that forming the divisorial test ideal $\widetilde{\tau}^{\mathrm{div}}((R,\Delta);\a^t)$ commutes with localization and completion. 
\begin{cor}[\textup{cf.\cite[Proposition 3.1, 3.2]{HT}}]\label{loc}\label{completion}
Let the notation be as in Lemma $\ref{key lemma}$.
\renewcommand{\labelenumi}{$(\arabic{enumi})$}
\begin{enumerate}
\item
Let $W$ be a multiplicatively closed subset of $R$, and let $\Delta_W$ and $\a_W$ be the images of $\Delta$ and $\a$ in $R_W$, respectively. Then
$$\widetilde{\tau}^{\mathrm{div}}((R_W,\Delta_W);\a_W^t) =\widetilde{\tau}^{\mathrm{div}}((R,\Delta);\a^t)R_W.$$
\item
Let $\widehat{R}$ be the $\m$-adic completion of $R$, and let $\widehat{\Delta}$ and $\widehat{\a}$ be the images of $\Delta$ and $\a$ in $\widehat{R}$, respectively.
Then
$$\widetilde{\tau}^{\mathrm{div}}((\widehat{R},\widehat{\Delta});\widehat{\a}^t)=\widetilde{\tau}^{\mathrm{div}}((R,\Delta);\a^t)\widehat{R}.$$
\end{enumerate}
\end{cor}

\begin{prop}\label{test}
Let $R$ be an F-finite normal domain of characteristic $p>0$ and $\Delta$ be a boundary on $\Spec R$. 
Let $\a \subseteq R$ be an ideal such that $\a \cap R^{\circ,\Delta} \ne \emptyset$ and $t>0$ be a real number. 
Let $E=\oplus_{\m} E_R(R/\m)$ be the direct sum, taken over all maximal ideals $\m$ of $R$, of the injective hulls of the residue fields $R/\m$. 
Fix an element $c \in \widetilde{\tau}^{\mathrm{div}}(R,\lfloor \Delta \rfloor) \cap R^{\circ,\Delta}$. 
Then for any $z \in E$, the following three conditions are equivalent to each other. 
\renewcommand{\labelenumi}{$(\arabic{enumi})$}
\begin{enumerate}
\item $z \in 0_{E}^{\mathrm{div}*(\Delta, \a^t)}$ 

\item $c\a^{\lceil tq \rceil}z^q=0$ in ${}^e\! R(\lceil q \Delta \rceil -\lfloor \Delta \rfloor) \otimes_R E$ 
for all $q=p^e$. 

\item There exist integers $n \geq 1$ and $n-1 \geq r \geq 0$ such that for each $q=p^e$, if $e \equiv r \mod n$, then $c\a^{\lceil tq \rceil}z^{q}=0$ in ${}^e\! R(\lceil q \Delta \rceil -\lfloor \Delta \rfloor)\otimes_R E$. 
\end{enumerate}
\end{prop}
\begin{proof}
We may assume that $R$ is local by the definition of $E$. For each $q=p^e$, we denote ${}^e\!R(\lceil q \Delta \rceil -\lfloor \Delta \rfloor)\otimes_R E$ by $\widetilde{\F}^{e,\Delta}(E)$. 
First we prove that the condition $(1)$ implies the condition $(2)$. 
Fix any $q=p^e$.  
Since $z \in 0_{E}^{\mathrm{div}*(\Delta, \a^t)}$, there exist $d \in R^{\circ,\Delta}$ and an integer $e_0 \geq 0$ such that  $d\a^{\lceil tq' \rceil}z^{q'}=0$ in $\widetilde{\F}^{e',\Delta}(E)$ for every $q'=p^{e'} \geq p^{e_0}$. 
By Lemma \ref{key lemma}, there exist an integer $e_1 \geq e_0$ and $R$-linear maps $\varphi^{(k)} \in \Hom_R(R((p^{k}-1)\lfloor \Delta \rfloor)^{1/p^{k}},R)$ for $e_1 \geq k \geq e_0$ such that $c = \sum_{k=e_0}^{e_1}\varphi^{(k)}(d^{1/p^{k}})$. 
The map $\varphi^{(k)}$ induces an ${}^e R$-linear map $\psi^{(k)}:\widetilde{\F}^{e+k,\Delta}(E) \to \widetilde{\F}^{e,\Delta}(E)$
sending $dz^{p^kq}$ to $c_{k}z^q$, where $c_{k}:=\varphi^{(k)}(d^{1/p^{k}})$ for $e_1 \geq k \geq e_0$. 
Since $d(\a^{\lceil tq \rceil})^{[p^k]} \subseteq  d\a^{\lceil tp^{k}q \rceil}$, $d(\a^{\lceil tq \rceil})^{[p^k]}z^{p^kq}=0$ in $\widetilde{\F}^{e+k,\Delta}(E)$ for every $e_1 \geq k \geq e_0$. 
Applying $\psi^{(k)}$ and summing up, we obtain
$$c\a^{\lceil tq \rceil}z^q=\sum_{k=e_0}^{e_1} c_{k}\a^{\lceil tq \rceil}z^q=\sum_{k=e_0}^{e_1}\psi^{(k)}(d(\a^{\lceil tq \rceil})^{[p^k]}z^{p^kq})=0$$
in $\widetilde{\F}^{e,\Delta}(E)={}^e\! R(\lceil q \Delta \rceil -\lfloor \Delta \rfloor) \otimes_R E$. 

The condition $(3)$ is the special case of $(2)$, so it remains to prove that the condition (3) implies the condition $(1)$.
\begin{claim}
There exists a positive integer $m$ such that for every $q=p^e$, one has an $R$-linear map $\phi_e:R((q-1)\lfloor \Delta \rfloor)^{1/q} \to R$ sending $1$ to $c^m$. 
\end{claim}
\begin{clproof}
By Corollary \ref{completion} and Proposition \ref{basic} (3), the pair $(R_c,\lfloor \Delta \rfloor_c)$ is divisorially F-regular. 
Then by \cite[Proposition 2.2]{HW}, there exist an integer $s \geq 1$ and an $R$-linear map $g:R((p-1)\lfloor \Delta \rfloor)^{1/p} \to R$ sending $1$ to $c^s$, because $R((p-1)\lfloor \Delta \rfloor)^{1/p}$ is a finite $R$-module.  

We prove by induction that for every $q=p^e$, there exists an $R$-linear map $\phi_e:R((q-1)\lfloor \Delta \rfloor)^{1/q} \to R$ sending $1$ to $c^{2s}$. 
When $e=1$, put $\phi_1=c^s \cdot g$. 
Suppose that the assertion holds for $q=p^e$. Then, by tensoring $\phi_e$ with $R((p-1)\lfloor \Delta \rfloor)$ and taking its $p$-th root, we have an $R^{1/p}$-module homomorphism $R((pq-1)\lfloor \Delta \rfloor)^{1/pq} \to R((p-1)\lfloor \Delta \rfloor)^{1/p}$ sending $1$ to $c^{2s/p}$. 
So, we compose this map with an $R^{1/p}$-linear map $R((p-1)\lfloor \Delta \rfloor)^{1/p} \to R((p-1)\lfloor \Delta \rfloor)^{1/p}$ sending $1$ to $c^{(p-2)s/p}$, and then with $g$:
$$R((pq-1)\lfloor \Delta \rfloor)^{1/pq} \to R((p-1)\lfloor \Delta \rfloor)^{1/p} \to R((p-1)\lfloor \Delta \rfloor)^{1/p} \xrightarrow{g} R$$
$$1 \mapsto c^{2s/p} \mapsto c^s \mapsto c^{2s}$$
This is the required map for $pq=p^{e+1}$. 
\end{clproof}

In order to prove the condition (1), it is enough to show that $c^{m+1} \a^{\lceil tq \rceil}z^q=0$ in $\widetilde{\F}^{e,\Delta}(E)$ for all $q=p^e$.
Fix any $q=p^e$ and take an integer $n \geq l  \geq 1$ such that $e+l \equiv r \mod n$. 
By the above claim, there exists an $R$-module homomorphism $\phi_l:R((p^l-1)\lfloor \Delta \rfloor)^{1/p^l} \to R$ sending $c$ to $c^{m+1}$. 
By tensoring $\phi_l$ with $R(\lceil q \Delta \rceil -\lfloor \Delta \rfloor)$ and taking its $q$-th root, 
we obtain an $R^{1/q}$-linear map $$R(\lceil p^lq \Delta \rceil -\lfloor \Delta \rfloor)^{1/p^lq} \to R(\lceil q \Delta \rceil -\lfloor \Delta \rfloor)^{1/q}.$$ This map induces an
${}^eR$-linear map $\psi_l: \widetilde{\F}^{e+l,\Delta}(E) \to \widetilde{\F}^{e,\Delta}(E)$
sending $c^{p^l}z^{p^lq}$ to $c^{m+1}z^q$. 
Since $c^{p^l}(\a^{\lceil tq \rceil})^{[p^l]} \subseteq c \a^{\lceil tp^lq \rceil}$, by the condition $(3)$, $c^{p^l}(\a^{\lceil tq \rceil})^{[p^l]}z^{p^lq}=0$ in  $\widetilde{\F}^{e+l,\Delta}(E)$.  
Applying $\psi_l$, we have
$$c^{m+1}\a^{\lceil tq \rceil}z^{q}=\psi_l(c^{p^l}(\a^{\lceil tq \rceil})^{[p^l]}z^{p^lq})=0 \in \widetilde{\F}^{e,\Delta}(E).$$ 
\end{proof}

\begin{rem}\label{test remark}
\renewcommand{\labelenumi}{(\roman{enumi})}
\begin{enumerate}
\item
We cannot replace $\widetilde{\F}^{e,\Delta}(E)={}^e\! R(\lceil q \Delta \rceil -\lfloor \Delta \rfloor)\otimes_R E$ by $\F^{e,\Delta}(E)={}^e\! R((q-1)\Delta) \otimes_R E$ in Proposition \ref{test}. 
\item
Let the notation be as in Proposition \ref{test}. 
In addition, let $I$ be an ideal of $R$. Then for any $x \in R$, by an argument similar to the proof of Proposition \ref{test}, the following three conditions are also equivalent to each other. 
\begin{enumerate}
\item $x \in I^{\mathrm{div}*(\Delta,\a^t)}$. 
\item $c\a^{\lceil tq \rceil}x^q \in I^{[q]}R(\lceil q \Delta \rceil -\lfloor \Delta \rfloor)$ for all $q=p^e$.  
\item There exist integers $n \geq 1$ and $n-1 \geq r \geq 0$ such that for each $q=p^e$,  
if $e \equiv r \mod n$, then $c\a^{\lceil tq \rceil}x^{q} \in I^{[q]}R(\lceil q \Delta \rceil -\lfloor \Delta \rfloor)$. 
\end{enumerate}
\end{enumerate}
\end{rem}

\begin{defn}\label{test element}
Let $R$ be an F-finite normal domain of characteristic $p>0$ and $D$ be a reduced divisor on $\Spec R$. 
Let $E=\oplus_{\m} E_R(R/\m)$ be the direct sum, taken over all maximal ideals $\m$ of $R$, of the injective hulls of the residue fields $R/\m$. 
We say that an element $c \in R^{\circ, D}$ is a \textit{divisorial $D$-test element} (for $E$) 
if $cz^q =0$ in $\F^{e,D}(E)$ for all $q=p^e$ whenever $z \in 0_{E}^{\mathrm{div}*D}$. 
\end{defn}

\begin{cor}\label{test cor}
Let $R$ and $D$ be as in Definition $\ref{test element}$.
\renewcommand{\labelenumi}{$(\arabic{enumi})$}
\begin{enumerate}
\item
$\widetilde{\tau}^{\mathrm{div}}(R,D) \cap R^{\circ,D}$ is equal to the set of divisorial $D$-test elements $($for $E)$. 

\item
Let $c \in R^{\circ,D}$ be an element such that the localized pair $(R_c,D_c)$ is divisorially F-regular. 
Then some power $c^n$ of $c$ is a divisorial $D$-test element $($for $E)$. 

\item
Let $\Delta$ be a boundary on $\Spec R$ and $\a \subseteq R$ be an ideal such that $\a \cap R^{\circ, \Delta} \ne \emptyset$. Then for any real number $t>0$, 
$\widetilde{\tau}^{\mathrm{div}}((R,\Delta);\a^t) \cap R^{\circ,\Delta} \ne \emptyset$. 

\item Assume that $D=\Div_R(f)$ is a reduced Cartier divisor on $\Spec R$. 
Let $\Delta$ be a boundary on $\Spec R$ whose support has no common components with the support of $D$, and let $\a \subseteq R$ be an ideal such that $\a \cap R^{\circ,\Delta+D} \ne \emptyset$. 
Assume in addition that $f$ lies in an ideal $\b \subseteq R$ such that $\b \cap R^{\circ,\Delta+D} \ne \emptyset$.  Then 
$$\widetilde{\tau}^{\mathrm{div}}((R,\Delta+D);\a^t) \subseteq \widetilde{\tau}^{\mathrm{div}}((R,\Delta);\a^t\b).$$

\end{enumerate}
\end{cor}

\begin{proof}
(1) It follows from Proposition \ref{test}  that  any element of $\widetilde{\tau}^{\mathrm{div}}(R,D) \cap R^{\circ,D}$ is a divisorial $D$-test element. 
Conversely, if $c$ is a  divisorial $D$-test element, then $cz^q=0$ in $\F^{e,D}(E)$ for all $z \in 0_E^{\mathrm{div}*D}$ and for all $q=p^e$. In particular, we have $c 0_E^{\mathrm{div}*D}=0$ in $E$, that is, $c \in\widetilde{\tau}^{\mathrm{div}}(R,D)$. 

(2) By Proposition \ref{basic} (3) and Corollary \ref{loc}, the localization $\widetilde{\tau}^{\mathrm{div}}(R,D)_c$ with respect to $c$ is trivial. 
This implies that some power $c^n$ of $c$ is in $\widetilde{\tau}^{\mathrm{div}}(R,D)$. 
Thus, $c^n$ is a divisorial $D$-test element by $(1)$. 

(3) Let $S$ be the ring corresponding to $\lfloor \Delta \rfloor$. 
Fix any element $c \in R^{\circ, \Delta}$ such that $R_c$ and $S_c$ are regular (one can choose such $c$, because $R$ is normal and $S$ is reduced). Then by \cite[Theorem 4.9]{HW}, $(R_c, \lfloor \Delta \rfloor_c)$ is divisorially F-regular. Thus, by (1) and (2), that some power $c^n$ of $c$ lies in $\widetilde{\tau}^{\mathrm{div}}(R,\lfloor \Delta \rfloor)$. The assertion follows from Proposition \ref{basic} (2). 

(4) Let $z \in 0_E^{\mathrm{div}*(\Delta,\a^t\b)}$, and fix any $c \in $$\widetilde{\tau}^{\mathrm{div}}(R,\Delta) \cap R^{\circ, \Delta+D}$ (one can choose such $c$ by (3)). 
By Proposition \ref{test}, $c\a^{\lceil tq \rceil}\b^qz^q=0$ in ${}^e\! R(\lceil q\Delta \rceil-\lfloor \Delta \rfloor) \otimes_R E$ for all $q=p^e$. Take any $d \in \b \cap R^{\circ,\Delta+D}$, and then by assumption, one has $cd \a^{\lceil tq \rceil}f^{q-1}z^q=0$ in ${}^e\! R(\lceil q\Delta \rceil-\lfloor \Delta \rfloor)  \otimes_R E$ for all $q=p^e$, which implies that $z \in 0_E^{\mathrm{div}*(\Delta+D,\a^t)}$. 
\end{proof}

\section{Geometric properties of divisorial test ideals}
\label{sec:4}
First we discuss a restriction property of $\widetilde{\tau}^{\mathrm{div}}((R,\Delta);\a^t)$ similar to Theorem \ref{restriction in char 0}.
To state our restriction theorem, we have to generalize the notion of the generalized test ideal $\widetilde{\tau}((R,\Delta);\a^t)$ to the case where $R$ is not necessarily normal.

Suppose that $R$ is a Noetherian equidimensional reduced ring satisfying $S_2$. 
Let $Q(R)$ be the total quotient ring of $R$ and $X:=\Spec R$. 
A fractional ideal $\I \subset Q(R)$ is called a \textit{Weil divisorial fractional ideal} of $R$ if $\I$ is principal in codimension one and reflexive. 
Since $R$ satisfies $S_2$, $R \subset Q(R)$ is a Weil divisorial fractional ideal. 
The product $\I \cdot \I' \subset Q(R)$ of Weil divisorial fractional ideals $\I,\I' \subset Q(R)$ is defined to be the reflexive hull of the product of the fractional ideals $\I\I' \subset Q(R)$. 
With this law, the set of Weil divisorial fractional ideals is a group which we denote by $\mathrm{WSh}(X)$. 
We write $\I^{(n)}$ for the product of $\I \in \mathrm{WSh}(X)$ with itself $n$-times in $\mathrm{WSh}(X)$. 
We say $\I \in \mathrm{WSh}(X)$ is \textit{effective} if $R \subseteq \I \subset Q(R)$. 
If $R$ has a canonical module $\omega_X$ (for example, $R$ is a homomorphic image of a Gorenstein local ring) and we fix an embedding $\omega_X \subset Q(R)$,  
then $\omega_X$ is a Weil divisorial fractional ideal precisely when $R$ is Gorenstein in codimension one (\cite[(16.3.3)]{Ko}).  
The reader is referred to \cite[Chapter 16]{Ko} for details. 

\begin{defn}
In the above situation, assume that $R$ is of characteristic $p>0$. 
Let $\I \in \mathrm{WSh}(X)$ be an effective Weil divisorial fractional ideal and $n$ be a positive integer. 
Let $\a \subseteq R$ be an ideal such that $\a \cap R^{\circ} \ne \emptyset$ and $t>0$ be a real number. 
Let $E=\oplus_{\m} E_R(R/\m)$ be the direct sum, taken over all maximal ideals $\m$ of $R$, of the injective hulls of the residue fields $R/\m$. 
Then the \textit{$(\I^{(1/n)},\a^t)$-tight closure} $0^{*(\I^{(1/n)},\a^t)}_E$
of the zero submodule in $E$ is defined to be the submodule of $E$ consisting of all elements $z \in E$ for which there exists $c \in R^{\circ}$ such that $c\a^{\lceil tq \rceil} \otimes z=0$ in ${}^e \I^{(\lfloor (q-1)/n \rfloor)} \otimes_R E$ for all large $q = p^e$,
where ${}^e \I^{(\lfloor (q-1)/n \rfloor)}$ is $\I^{(\lfloor (q-1)/n \rfloor)}$ itself but viewed as an $R$-module via the $e$-times Frobenius map $F^e:R \to \I^{(\lfloor (q-1)/n \rfloor)}$. 
We define the ideal $\widetilde{\tau}((R,\I^{(1/n)});\a^t)$ by 
$$\widetilde{\tau}((R,\I^{(1/n)});\a^t)=\mathrm{Ann}_R(0^{*(\I^{(1/n)},\a^t)}_E) \subseteq R. $$ 
\end{defn}

\begin{propdef}[\textup{\cite[16.5 Proposition-Definition]{Ko}}]\label{different}
Let $R$ be a normal domain and $B$ be an effective $\Q$-divisor on $X:=\Spec R$. 
Let $i:D \hookrightarrow X$ be a reduced subscheme of pure codimension one which is Gorenstein in codimension one and satisfies $S_2$. 
Let $S$ be the ring corresponding to $D$,  and suppose that $R$ $($resp. $S)$ has a canonical module $\omega_X \subset Q(R)$ $($resp. $\omega_D \subset Q(S))$.
Assume that $K_X+D+B$ is $\Q$-Cartier, and let $r \geq 1$ be the least integer such that $r(K_X+D+B)$ is Cartier.  
Assume in addition that 
$$\mathrm{Supp}(\omega_X(D+B)):=\{x \in X \mid \omega_X(D+B)_x \ne \sO_{X,x}\}$$ has no common components with the support of $D$. 
Then there exists a naturally defined effective Weil divisorial fractional ideal $\I_B \in \mathrm{WSh}(D)$ so that:
$$\omega_D^{(r)} \cdot \I_B =i^* \omega_X^{(r)}(r(D+B)).$$
We use the formal exponential notation $\I_B^{(1/r)}$ for the pair $(r,  \I_B)$ and call $\I_B^{(1/r)}$ as the \textit{different} of $B$ on $D$. 
\end{propdef}

\begin{rem}
If $R$ is $\Q$-Gorenstein and Cohen-Macaulay, $D$ is a Cartier reduced divisor on $X$ and $B=0$, then $\I_B$ is equal to $S$ in Definition \ref{different}. The reader is referred to \cite[16.4 Proposition]{Ko} for details. 
\end{rem}

Now we state our restriction theorem for the ideal $\widetilde{\tau}^{\mathrm{div}}((R,\Delta);\a^t)$ (compare this with Theorem \ref{restriction in char 0}). 
\begin{thm}\label{restriction}
Let $(R,\m)$ be a $d$-dimensional F-finite normal local ring of characteristic $p>0$ which is a homomorphic image of a Gorenstein ring. 
Let $\Delta$ be a boundary on $X:=\Spec R$ such that $K_X+\Delta$ is $\Q$-Cartier. 
Let $r \geq 1$ be the least integer such that $r(K_X+\Delta)$ is Cartier, and assume that $r$ is not divisible by $p$. 
Put $D=\lfloor \Delta \rfloor$ and $B=\{\Delta\}$, and assume in addition that $D$ is Gorenstein in codimension one and satisfies $S_2$. 
Let $S$ be the local ring corresponding to $D$ and $\mathcal{I}_B^{(1/r)}$ be the different of $B$ on $D$. 
Then for any real number $t>0$ and for any ideal $\a \subseteq R$ such that $\a \cap R^{\circ, \Delta} \ne \emptyset$, one has 
$$\widetilde{\tau}((S,\I_B^{(1/r)});{\a S}^t) = \widetilde{\tau}^{\mathrm{div}}((R,\Delta); \a^t)S.$$
In particular, if $\widetilde{\tau}^{\mathrm{div}}(R,\Delta)=R$, then $S$ is strongly F-regular $($hence normal and Cohen-Macaulay$)$. 
\end{thm}

\begin{proof}
We apply a strategy similar to \cite[Theorem 4.9]{HW}. 
First note that the ideal $\widetilde{\tau}((S,\I_B^{(1/r)});{\a S}^t)$ commutes with completion by an argument similar to that of \cite[Proposition 3.2]{HT}. 
So, by virtue of Lemma \ref{completion}, 
we may assume without loss of generality that $R$ is complete. 
Let $E_R=E_R(R/\m)$ and $E_S=E_{S}({S}/\m S)$ be the injective hulls of the residue fields of $R$ and $S$, respectively. 
Since $R$ and $S$ satisfy $S_2$,  by \cite{A1} and \cite{A2}, $E_R \cong H^d_{\m}(\omega_X)$ and $E_S  \cong H^{d-1}_{\m S}(\omega_{D})$. 
Also, we can view $E_{S}$ as a submodule of $E_R$ via the isomorphism $E_{S} \cong (0:R(-D))_{E_R} \subset E_R$. 
\begin{claim}
$$0^{\mathrm{div}*(\Delta,\a^t)}_{E_R} \cap E_{S} = 0^{*(\I_B^{(1/r)}; {\a S}^t)}_{E_{S}}.$$
\end{claim}
\begin{clproof}
Fix any $d \in R(-\lceil B \rceil) \cap R^{\circ, \Delta}$. 
One can choose a divisorial $D$-test element $c \in  R^{\circ,\Delta}$ whose image is in $\widetilde{\tau}(S) \cap S^{\circ}$. 
We stick to the powers $q=p^e$ such that $q-1$ is divisible by $r$. 
Note that there exists an integer $n$ such that if a positive integer $e$ is a multiple of $n$, then $p^e-1$ is divisible by $r$ (because  $p$ does not divide $r$). 
Then $(q-1)(K_X+\Delta)$ is a Cartier divisor on $X$ and $\I_B^{((q-1)/r)}$ is an effective Weil divisorial fractional ideal of $S$. 
The inclusion map $E_S \hookrightarrow E_R$ induces the following commutative diagram of local cohomology modules, where the vertical maps are the induced $e$-times Frobenius maps:
$$\xymatrix{
0 \ar[r] &H^{d-1}_{\m S}(\omega_{D}) \ar[d]^{\overline{cd}F_S^e} \ar[r]^{} & H^d_{\m}(\omega_X) \ar[d]^{cdF_R^e}\\
0 \ar[r] & H^{d-1}_{\m S}(\omega_{D}^{(q)} \cdot \I_B^{((q-1)/r)}) \ar[r]^{} & H^d_{\m}(\omega_X^{(q)}((q-1)\Delta))\\
}$$

Let $\xi \in 0^{\mathrm{div}*(\Delta,\a^t)}_{E_R} \cap E_{S}$. 
Since $c$ is a divisorial $D$-test element, by Proposition \ref{test} and Corollary \ref{test cor}, 
$c \a^{\lceil tq \rceil}F^e_R(\xi)=0$ in $H^d_{\m}(\omega_X^{(q)}(\lceil q\Delta \rceil-D))$.  
Applying the multiplication map by $d$
$$H^d_{\m}(\omega_X^{(q)}(\lceil q\Delta \rceil-D)) \xrightarrow{\times d} H^d_{\m}(\omega_X^{(q)}((q-1)\Delta)),$$ 
one has 
$cd \a^{\lceil tq \rceil}F^e_R(\xi)=0$ in $H^d_{\m}(\omega_X^{(q)}((q-1)\Delta))$.  
Then, by the above commutative diagram, $\overline{cd}\a^{\lceil tq \rceil}F^e_{S}(\xi)=0$ in $H^{d-1}_{\m S}(\omega_{D}^{(q)} \cdot \I_B^{((q-1)/r)})$ for all $q=p^e$ such that $e$ is a multiple of $n$. 
Since $\overline{cd} \in \widetilde{\tau}(S) \cap S^{\circ}$, it follows from an argument similar to the proof of Proposition \ref{test} $(3) \Rightarrow (1)$ that $\xi$ lies in $0^{*(\I_B^{(1/r)}, {\a S}^t)}_{E_{S}}$. 

Consequently, we have $0^{\mathrm{div}*(\Delta,\a^t)}_{E_R} \cap E_{S} \subseteq 0^{*(\I_B^{(1/r)}, {\a S}^t)}_{E_{S}}$. 
The converse argument just reverses this.  
\end{clproof}
We continue the proof of Theorem \ref{restriction}. 
Since $R$ is complete,  one has $0^{\mathrm{div}*(\Delta,\a^t)}_{E_R}=(0:\widetilde{\tau}^{\mathrm{div}}((R,\Delta);\a^t))_{E_R}$. 
Hence 
\begin{align*}
0^{\mathrm{div}*(\Delta,\a^t)}_{E_R}\cap E_S &=(0:(\widetilde{\tau}^{\mathrm{div}}((R,\Delta);\a^t)+R(-D)))_{E_R}\\
&=\left( 0: \frac{\widetilde{\tau}^{\mathrm{div}}((R,\Delta);\a^t)+R(-D)}{R(-D)} \right)_{E_S}\\
&=(0:\widetilde{\tau}^{\mathrm{div}}((R,\Delta);\a^t)S)_{E_S}.
\end{align*}
By the above claim, we conclude that 
$$\widetilde{\tau}((S,\I_B^{(1/r)});{\a S}^t) =\mathrm{Ann}_{S}  (0^{\mathrm{div}*(\Delta,\a^t)}_{E_R}\cap E_S)=\widetilde{\tau}^{\mathrm{div}}((R,\Delta); \a^t)S.$$
\end{proof}

Next we discuss subadditivity properties, analogous to Theorem \ref{sub in char 0}, of the ideals $\widetilde{\tau}((R,\Delta);\a^t)$ and $\widetilde{\tau}^{\mathrm{div}}((R,\Delta);\a^t)$, respectively. 
\begin{thm}\label{subadditivity}
\renewcommand{\labelenumi}{$(\arabic{enumi})$}
\begin{enumerate}
\item
Let $(R,\m)$ be a $d$-dimensional F-finite regular local ring of characteristic $p>0$. 
Then  
$$\widetilde{\tau}(R,\a^s\b^t) \subseteq \sum_{\lambda+\mu=d} \widetilde{\tau}(R,\a^s\m^{\lambda})\widetilde{\tau}(R,\b^t\m^{\mu})$$
for any nonzero ideals $\a, \b \subseteq R$ and for any real numbers $s,t>0$. 
\item
Let $R$ be an F-finite regular domain of characteristic $p>0$ or an affine normal domain over a perfect field $K$ of positive characteristic $p$. Let $\mathfrak{J}$ be the Jacobian ideal of $R$ over $K$  $(\mathfrak{J}=R$ if $R$ is regular$)$.  
Let $\Delta_1, \Delta_2$ be boundaries on $\Spec R$ such that $\Delta_1+\Delta_2$ is also a boundary. 
Furthermore, assume that there exists an integer $r \geq 1$ not divisible by $p$ such that $r\Delta_1$ or $r\Delta_2$ is Cartier. 
Then 
$$\mathfrak{J} \widetilde{\tau}^{\mathrm{div}}((R,\Delta_1+\Delta_2);\a^s\b^t) 
\subseteq  \widetilde{\tau}^{\mathrm{div}}((R,\Delta_1);\a^s)\widetilde{\tau}^{\mathrm{div}}((R,\Delta_2);\b^t)$$
for any ideals $\a, \b \subseteq R$ such that $\a\b \cap R^{\circ,\Delta_1+\Delta_2} \ne \emptyset$ and for any real numbers $s,t>0$. 
\end{enumerate}
\end{thm}

\begin{proof}
$(1)$ We employ the same strategy as that of the proof of \cite[Theorem 6.10]{HY} which can be traced back to the method in \cite{DEL}. 
By virtue of Corollary \ref{completion}, we may assume that $R$ is a complete regular local ring of positive characteristic $p$, that is, $R=k[[x_1, \dots, x_d]]$ where $k$ is a field of characteristic $p$. 
Let $T= R \widehat{\otimes}_k R=k[[x_1, \dots, x_d, y_1, \dots, y_d]]$ and $\Delta:T \to R$ be the diagonal map. 
$\Ker(\Delta)$ is generated by $x_1-y_1, \dots, x_d-y_d$, and put $D_i=\mathrm{div}_T(x_i-y_i)$ for all $i=1,\dots,d$. 
Also, let $p_1$ (resp. $p_2$) $:\Spec T \to \Spec R$ be the first (resp. second) projection. 
Applying Theorem \ref{restriction} and Corollary \ref{test cor} (4) to the diagonal map $\Delta:T \to R$, 
one has the following inclusion: 
\begin{align*}
\widetilde{\tau}(R,\a^s\b^t) &=\Delta(\widetilde{\tau}^{\mathrm{div}}((T,D_1+\dots+D_d);(p_1^{-1}\a)^s (p_2^{-1}\b)^t))\\
& \subseteq \Delta(\widetilde{\tau}(T,(p_1^{-1}\a)^s (p_2^{-1}\b)^t\m_T^d)).
\end{align*}
Here, note that $\m_T=p_1^{-1}\m+p_2^{-1}\m$. 
It therefore follows from \cite[Theorem 3.1]{Ta3} and \cite[Proposition 4.4]{HY} that 
\begin{align*}
\widetilde{\tau}(T,(p_1^{-1}\a)^s (p_2^{-1}\b)^t\m_T^d) & =
\sum_{\lambda+\mu=d}\widetilde{\tau}(T,(p_1^{-1}\a)^s(p_1^{-1}\m)^{\lambda} (p_2^{-1}\b)^t(p_2^{-1}\m)^{\mu})\\
& \subseteq \sum_{\lambda+\mu=d}p_1^{-1}\widetilde{\tau}(R,\a^s\m^{\lambda}) p_2^{-1} \widetilde{\tau}(R,\b^t\m^{\mu}).
\end{align*}
Thus, we have the required inclusion 
$$\widetilde{\tau}(R,\a^s\b^t) \subseteq \sum_{\lambda+\mu=d} \widetilde{\tau}(R,\a^s\m^{\lambda})\widetilde{\tau}(R,\b^t\m^{\mu}).$$

$(2)$ 
We employ the same strategy as that of the proof of \cite[Theorem 2.7]{Ta3}. 
By assumption, we may assume that there exists an integer $r \geq 1$ not divisible by $p$ such that $r\Delta_1$ is Cartier. 
By virtue of Corollary \ref{completion}, we may also assume that $(R,\m)$ is a complete local ring. 
Let $E=E_R(R/\m)$ be the injective hull of the residue field of $R$, and we denote ${}^e\!R(\lceil q (\Delta_1+\Delta_2)  \rceil-\lfloor \Delta_1+\Delta_2 \rfloor) \otimes_R E$ (resp. ${}^e\!R(\lceil q \Delta_i  \rceil-\lfloor \Delta_i \rfloor) \otimes_R E$) by $\widetilde{\F}^{e,\Delta_1+\Delta_2}(E)$ (resp.  $\widetilde{\F}^{e,\Delta_i}(E)$) for each $q=p^e$.  

\begin{claim} For any ideal $I \subseteq R$, one has 
$$I^{\mathrm{div}*(\Delta_1,\a^s)} \widetilde{\tau}^{\mathrm{div}}((R,\Delta_1+\Delta_2);\a^s\b^t)
\subseteq I \widetilde{\tau}^{\mathrm{div}}(R,\Delta_2);\b^t).$$
\end{claim}
\begin{clproof}
By an argument similar to the proof of \cite[Proposition 2.2]{Ta3}, it is enough to show that
$$(0_E^{\mathrm{div}*(\Delta_1+\Delta_2,\a^s\b^t)}:I^{\mathrm{div}*(\Delta_1,\a^s)})_E \supseteq (0_E^{\mathrm{div}*(\Delta_2,\b^t)}:I)_E.$$
Let $z \in (0_E^{\mathrm{div}*(\Delta_2,\b^t)}:I)_E$, and take a divisorial $\lfloor \Delta_1+\Delta_2 \rfloor$-test element  $c \in R^{\circ,\Delta_1+\Delta_2}$.  
Since $c$ is also a divisorial $\lfloor \Delta_2 \rfloor$-test element,  by Proposition \ref{test} and Corollary \ref{test cor}, one has $c\b^{\lceil tq \rceil}I^{[q]}z^q=0$ in $\widetilde{\F}^{e,\Delta_2}(E)$ for all $q=p^e$. 
Then it follows that $c\b^{\lceil tq \rceil}I^{[q]}R(\lceil q \Delta_1 \rceil- \lfloor \Delta_1 \rfloor)z^q=0$ in $\widehat{\F}^{e,\Delta_1+\Delta_2}(E)$. 
On the other hand, since $c$ is a divisorial $\lfloor \Delta_1 \rfloor$-test element, by Remark \ref{test remark} (ii) and Corollary \ref{test cor}, one has $c\a^{\lceil sq \rceil}(I^{\mathrm{div}*(\Delta_1,\a^s)})^{[q]} \subseteq I^{[q]}R(\lceil q \Delta_1 \rceil- \lfloor \Delta_1 \rfloor)$ for all $q=p^e$.
Thus, 
$$c^2\a^{\lceil sq \rceil}\b^{\lceil tq \rceil}(I^{\mathrm{div}*(\Delta_1,\a^s)})^{[q]}z^q 
=0 \in \widetilde{\F}^{e,\Delta_1+\Delta_2}(E)$$
for all $q=p^e$, that is, $ z \in (0_E^{\mathrm{div}*(\Delta_1+\Delta_2,\a^s\b^t)}:I^{\mathrm{div}*(\Delta_1,\a^s)})_E$. 
\end{clproof}

Thanks to the above claim, it suffices to show that the Jacobian ideal $\mathfrak{J}$ is contained in  $\widetilde{\tau}^{\mathrm{div}}((R,\Delta_1);\a^s)^{\mathrm{div}*(\Delta_1,\a^s)}$. 
Let $c \in \mathfrak{J}$, and fix any $q=p^e$ such that $(q-1)\Delta_1$ is Cartier. 
Note that there exists an integer $n \geq 1$ such that if a positive integer $e$ is a multiple of $n$, then $(p^e-1)\Delta_1$ is Cartier, because $p$ does not divide $r$. 
As in the proof of \cite[Lemma 2.6]{Ta3},  we may assume that there exists a finite free $R$-algebra $S_q$ such that $R \subseteq S_q \subseteq R^{1/q}$ and $cR^{1/q} \subseteq S_q$. 
Let $d \in R^{\circ, \Delta_1}$ be a divisorial $\lfloor \Delta_1 \rfloor$-test element, and then by Proposition \ref{test} and Corollary \ref{test cor}, $d \a^{\lceil sq \rceil}z^q=0$ in $\widetilde{\F}^{e,\Delta_1}(E)$ for all $z \in 0_E^{\mathrm{div}*(\Delta_1,\a^s)}$.  
Applying the multiplication map $\widetilde{\F}^{e,\Delta_1}(E) \xrightarrow{\times d'} \F^{e,\Delta_1}(E)$
with  $d' \in R(- \{\Delta_1\}) \cap R^{\circ, \Delta_1}$,  
one has 
$dd' \a^{\lceil sq \rceil}z^q=0$ in $\F^{e,\Delta_1}(E)$.  
Let $F^e_E:E \to R^{1/q} \otimes_R E$ be the $e$-times iterated Frobenius map induced on $E$, and 
write down $(q-1)\Delta_1=\Div_R(f_q)$ with $f_q \in R$. 
Since $0_E^{\mathrm{div}*(\Delta_1,\a^s)}=(0:\widetilde{\tau}^{\mathrm{div}}((R,\Delta_1);\a^s))_E$, 
 it follows that $(dd')^{1/q}f_q^{1/q}\a^{\lceil sq \rceil/q}F^e_E((0:\widetilde{\tau}^{\mathrm{div}}((R,\Delta_1);\a^s))_E)=0$ in $R^{1/q} \otimes_R E$.  
Here we consider the following $R$-module homomorphism:
$$\phi_q: E \xrightarrow{F^e_{E}}  R^{1/q} \otimes_R E \xrightarrow{\times c} S_q \otimes_R E.$$
Then $(dd')^{1/q}f_q^{1/q}\a^{\lceil sq \rceil/q}\phi_q((0:\widetilde{\tau}^{\mathrm{div}}((R,\Delta_1);\a^s))_E)=0$ in $S_q \otimes_R E$, because $\phi_q$ factors through $F^e_{E}$.  
On the other hand, as an $R$-module, 
$\phi_q((0:\widetilde{\tau}^{\mathrm{div}}((R,\Delta_1);\a^s))_E)$ generates $c (0 : \widetilde{\tau}^{\mathrm{div}}((R,\Delta_1);\a^s)S_q)_{S_q \otimes_R E}$ in $S_q \otimes_R E$, because $\phi_q$ is factorized into $E \to S_q \otimes_R E \xrightarrow{\times c} S_q \otimes_R E$ and $S_q$ is flat over $R$. 
Thus, 
$$(dd')^{1/q}f_q^{1/q}\a^{\lceil sq \rceil/q} c (0 : \widetilde{\tau}^{\mathrm{div}}((R,\Delta_1);\a^s)S_q)_{S_q \otimes_R E}=0$$ in $S_q \otimes_R E$.  
Since $S_q$ is a finitely generated free $R$-module,  we have an isomorphism $\Hom_{R}(S_q, E) \cong S_q \otimes_R E$.  
We apply Matlis duality to $S_q$ via this isomorphism so that 
\begin{align*}
(dd')^{1/q}f_q^{1/q}\a^{\lceil sq \rceil/q} c & \subseteq \Ann_{S_q}( (0 : \widetilde{\tau}^{\mathrm{div}}((R,\Delta_1);\a^s)S_q)_{S_q \otimes_R E})\\
&=\widetilde{\tau}^{\mathrm{div}}((R,\Delta_1);\a^s)S_q\\
& \subseteq \widetilde{\tau}^{\mathrm{div}}((R,\Delta_1);\a^s)R^{1/q}.
\end{align*}
This implies that $dd'\a^{\lceil sq \rceil} c^q \subseteq \widetilde{\tau}^{\mathrm{div}}((R,\Delta_1);\a^s)^{[q]}R((q-1)\Delta_1)$ for all $q=p^e$ such that $e$ is a multiple of $n$. 
By Remark \ref{test remark} (ii), we conclude that $c$ belongs to $\widetilde{\tau}^{\mathrm{div}}((R,\Delta_1);\a^s)^{\mathrm{div}*(\Delta_1,\a^s)}$.
\end{proof}

\begin{eg}\label{t.c. example}
\renewcommand{\labelenumi}{(\roman{enumi})}
\begin{enumerate}
\item
 Let $R=k[[x_1, \dots, x_d]]$ be a $d$-dimensional complete regular local ring of characteristic $p>0$ and $\Delta=\Div_R(x_1 \dots x_d)$. 
Then we can easily check that $\widetilde{\tau}^{\mathrm{div}}(R,\Delta)=(x_1, \dots, x_d)$. 
\item
Let $R=k[[x_1, \dots, x_d]]$ be a complete regular local ring of characteristic $p>0$ with $d \geq 2$ and $\Delta=\Div_R(x_1^n+ \cdots +x_d^n)$. 
Assume that the integer $n$ is not divisible by $p$ and $p \geq (n-1)(d-2)$. 
Then, by \cite[Theorem 6.4]{Hu} and Theorem \ref{restriction}, one has that $\widetilde{\tau}^{\mathrm{div}}(R,\Delta)=(x_1, \dots, x_d)^{n-d+1}$. 

\item 
Let $R$ be a two-dimensional complete $\Q$-Gorenstein normal local ring of characteristic $p>0$ and $\Delta$ be a reduced Cartier divisor on $\Spec R$.  
Since the test ideal of a one-dimensional reduced local ring is equal to its conductor ideal, by Theorem \ref{restriction}, $\widetilde{\tau}^{\mathrm{div}}(R,\Delta)$ is the lift of the conductor ideal $\mathfrak{c}(\Delta)$ of $\Delta$. 

\item 
Let $(R,\m)$ be a $d$-dimensional Cohen-Macaulay $\Q$-Gorenstein normal local ring essentially of finite type over a perfect field $K$ of characteristic $p>0$, and let  $\Delta$ be a reduced Cartier divisor on $\Spec R$. Let $S$ be the local ring corresponding to $\Delta$ and $\mathfrak{J}(S/K) \subset S$ be the Jacobian ideal of $S$ over $K$. 
Then by arguments similar to the proof of \cite[Theorem 4.6]{Ta3} and Theorem \ref{restriction}, $$\mathfrak{J}(S/K) \subseteq \widetilde{\tau}^{\mathrm{div}}((R,\Delta);\m^{d(1-\epsilon)})S$$ for all $1 \gg \epsilon >0$. 
In particular, $\widetilde{\tau}^{\mathrm{div}}(R,\Delta)$ contains the lift of the Jacobian ideal $\mathfrak{J}(S/K)$. 
\end{enumerate}
\end{eg}

\section{An interpretation of adjoint ideals via divisorial tight closure}
\label{sec:5}
In this section, we prove the correspondence between the divisorial test ideal $\widetilde{\tau}^{\mathrm{div}}((R,\Delta);\a^t)$ and the adjoint ideal $\adj((X,\Delta);\a^t)$.  
First of all, we show that $\widetilde{\tau}^{\mathrm{div}}((R,\Delta);\a^t)$ is contained in $\adj((X,\Delta);\a^t)$ in fixed prime characteristic (if a log resolution of $((X,\Delta);\a)$ exists).  
\begin{thm}\label{I < adj}
Let $(R,\m)$ be an F-finite normal local ring of characteristic $p>0$ which is a homomorphic image of a Gorenstein ring, and let $\Delta$ be a boundary $($i.e., $\Delta=\sum_i d_i \Delta$ is an $\R$-divisor with $0 \leq d_i \leq 1)$ on $X:=\Spec R$ such that $K_X+\Delta$ is $\R$-Cartier. 
Let $\a \subseteq R$ be an ideal such that no component of $\lfloor \Delta \rfloor$ is contained in the zero locus of $\a$, and let $t >0$ be a fixed real number.   
If $f:\widetilde{X} \to X$ is a proper birational morphism from a normal scheme $\widetilde{X}$ such that $\a \sO_{\widetilde{X}}=\sO_{\widetilde{X}}(-Z)$ is invertible, then one has an inclusion 
$$\widetilde{\tau}^{\mathrm{div}}((R,\Delta);\a^t) \subseteq H^0(\widetilde{X}, \sO_{\widetilde{X}}( K_{\widetilde{X}}- \lfloor f^*(K_X+\Delta)+tZ \rfloor +f^{-1}_*\lfloor \Delta \rfloor)).$$
\end{thm}
\begin{proof}
The proof is almost the same as those of \cite[Theorem 3.3]{HW} and \cite[Theorem 2.13]{Ta}.
\end{proof}

\begin{lem}[\textup{\cite[Proposition 3.6, Corollary 3.8]{Ha1}}]\label{char p>>0}
Let $(R,\m)$ be a $d$-dimensional normal local ring of essentially of finite type over a perfect field $\kappa$ of characteristic $p>0$. 
Let $f:\widetilde{X} \to X:=\Spec R$ be a log resolution of $X$ and $E$ be an $f$-ample $\R$-divisor on $\widetilde{X}$ whose fractional part has simple normal crossing support. 
Denote the closed fiber of $f$ by $Z=f^{-1}(\{m\})$, and assume that  $(R,\m)$ is the localization at any prime ideal of a finitely generated $\kappa$-algebra which is a reduction modulo $p \gg 0$ as well as $f:\widetilde{X} \to \Spec R$, $E$ and $Z$. 
Then the $e$-times Frobenius map 
$$F^e:H^d_Z(\widetilde{X}, \sO_{\widetilde{X}}(-E)) \to H^d_Z(\widetilde{X}, \sO_{\widetilde{X}}(-qE))$$
is injective for all $q=p^e$. 
\end{lem}

We prove that the adjoint ideal $\adj((X,\Delta);\a^t)$ coincides, after reduction to characteristic $p \gg 0$, with the divisorial test ideal $\widetilde{\tau}^{\mathrm{div}}((R,\Delta);\a^t)$, by making use of the above lemma. 
\begin{thm}
Let $(R,\m)$ be a $d$-dimensional normal local ring essentially of finite type over a perfect field of prime characteristic $p$, and let $\Delta$ be a boundary $($i.e., $\Delta=\sum_i d_i \Delta_i$ is an $\R$-divisor with $0 \leq d_i \leq 1)$ on $X:=\Spec R$ such that $K_X+\Delta$ is $\R$-Cartier. 
Let $\a \subseteq R$ be an ideal such that no component of $\lfloor \Delta \rfloor$ is contained in the zero locus of $\a$, and let $t >0$ be a fixed real number. 
Assume that $((R,\Delta);\a)$ is reduced from characteristic zero to characteristic $p \gg 0$, together with a log resolution $f:\widetilde{X} \to X$ of $((X,\Delta);\a)$ giving the adjoint ideal $\adj((X,\Delta);\a^t)$ $($see Definition \ref{adjoint} for the definition of $\adj((X,\Delta);\a^t))$. Then
$$\adj((X,\Delta);\a^t)=\widetilde{\tau}^{\mathrm{div}}((R,\Delta);\a^t).$$
\end{thm}
\begin{proof}
The proof is an improvement of the proof of \cite[Theorem 3.2]{Ta}\footnote{The proof of  \cite[Theorem 3.2]{Ta} has a small gap, but the result itself is valid. It follows from almost the same argument as the proof of Theorem 5.3.}. 
By Theorem \ref{I < adj}, it is enough to prove  an inclusion $\adj((X,\Delta);\a^t) \subseteq \widetilde{\tau}^{\mathrm{div}}((R,\Delta);\a^t)$.  
Since $\adj((X,\Delta);\a^t)$ does not change after an arbitrarily small perturbation of $\Delta$, we may assume that $K_X+\Delta$ is $\Q$-Cartier. We may also assume that $t$ is a rational number. 

First we consider the situation in characteristic zero before reducing to characteristic $p \gg 0$. 
In characteristic zero, we fix any $c \in R^{\circ, \Delta}$ lying in $\adj(X,\Delta) \cap R(-\{\Delta\}) \cap \a$ such that the localization $R_c$ is regular  (one can choose such $c$, because $R$ is normal). 
Take a log resolution $f:\widetilde{X} \to X$ of $((X,\Delta);c\a)$ so that $f^{-1}_*\lfloor \Delta \rfloor$ is smooth and $\a\sO_{\widetilde{X}}=\sO_{\widetilde{X}}(-F)$ for an effective integral divisor $F$ on $\widetilde{X}$. 
One can choose an $f$-ample $\Q$-divisor $H$ on $\widetilde{X}$ with simple normal crossing support so that $tF-H$ is an effective $\Q$-divisor on $\widetilde{X}$ whose support has no common components with the support of $\lfloor \Delta \rfloor$ and in addition
$$\lfloor f^*(K_X+\Delta) +tF \rfloor = \lfloor f^*(K_X+\Delta)+\epsilon \ \Div_{\widetilde{X}}(c)+tF-H \rfloor.$$ 
for a sufficiently small rational number $1 \gg \epsilon>0$. 
We denote $D=\lfloor \Delta \rfloor$, $\widetilde{D}=f^{-1}_*\lfloor \Delta \rfloor$ and $A= f^*(K_X+\Delta)+\epsilon \ \Div_{\widetilde{X}}(c)+tF-H$. 

\begin{cln}
For each integer $m \geq 1$, a natural exact sequence  
$$0 \to \sO_{\widetilde{X}}(mA-\widetilde{D}) \to \sO_{\widetilde{X}}(mA) \to \sO_{\widetilde{D}}(mA) \to 0$$
induces the exact sequence of local cohomology modules
$$0 \to H^{d-1}_{\widetilde{D} \cap Z}(\widetilde{D},\sO_{\widetilde{D}}(m A)) \to H^d_Z(\widetilde{X}, \sO_{\widetilde{X}}(mA-\widetilde{D})) \to H^{d}_Z(\widetilde{X}, \sO_{\widetilde{X}}(mA)) \to 0,$$
where $Z=f^{-1}(\{\m \})$ is the closed fiber of $f$. 
\end{cln}
\begin{clnproof}
It suffices to show that $H^{d-1}_Z(\widetilde{X}, \sO_{\widetilde{X}}(mA))=0$ for all $m \geq 1$. 
It immediately follows from Kawamata-Viehweg vanishing theorem \cite[Theorem 1-2-3]{KMM},  because $H^{d-1}_Z(\widetilde{X}, \sO_{\widetilde{X}}(mA))$ is the Matlis dual of $H^1(\widetilde{X}, \sO_{\widetilde{X}}(K_{\widetilde{X}}+\lceil -mA \rceil))$ and $-mA$ is an $f$-ample $\Q$-divisor whose fractional part has  simple normal crossing support.  
\end{clnproof}

Now we reduce the entire setup as above to characteristic $p \gg 0$ and switch the notation to denote things after reduction modulo $p$.  
Assume that $\a$ is generated by $r$ elements. 
Also, since $H-tF$ is an $f$-ample $\Q$-divisor on $\widetilde{X}$, putting $G:=K_{\widetilde{X}}-f^*(K_X+\Delta)+\widetilde{D}$, we see that $\mathcal{M}:=\bigoplus_{n \geq 0} H^0(\widetilde{X},\sO_{\widetilde{X}}(\lceil G-ntF+nH \rceil))$ is a finitely generated module over $\mathcal{R}:=\bigoplus_{n \geq 0} H^0(\widetilde{X},\sO_{\widetilde{X}}(n(H-tF)))$. 
So suppose that $\mathcal{M}$ is generated in degree $\leq n_0$. 
Since $R_c$ is regular and $(\Spec R_c, D_c)$ is plt, by \cite[Theorem 5.50]{KM} and \cite[Theorem 4.8]{HW}, $(R_c, D_c)$ is divisorially F-regular. 
It therefore follows from Corollary \ref{test cor} that some power $c^s$ is a divisorial $D$-test element. 
For all $q=p^e$, we consider the following commutative diagram with exact rows, where the vertical maps are induced by the $e$-times Frobenius map.  
$$\xymatrix{
0 \ar[r] & \sO_{\widetilde{X}}(A-\widetilde{D}) \ar[d]^{F_{\widetilde{X}}^e} \ar[r] & \sO_{\widetilde{X}}(A) \ar[d]^{F_{\widetilde{X}}^e} \ar[r] & \sO_{\widetilde{D}}(A) \ar[d]^{F_{\widetilde{D}}^e} \ar[r] & 0 \\
0 \ar[r] & \sO_{\widetilde{X}}(qA-\widetilde{D}) \ar[r] & \sO_{\widetilde{X}}(qA) \ar[r] & \sO_{\widetilde{D}}(qA) \ar[r] & 0 \\
}$$
By Claim $1$, this diagram induces the following commutative diagram of local cohomology modules:
\begin{small}
$$\xymatrix{
0 \ar[r] & H^{d-1}_{\widetilde{D} \cap Z}(\widetilde{D},\sO_{\widetilde{D}}(A)) \ar[d]^{F_{\widetilde{D}}^e} \ar[r] & H^d_Z(\widetilde{X}, \sO_{\widetilde{X}}(A-\widetilde{D})) \ar[d]^{F_{\widetilde{X}}^e} \ar[r] & H^{d}_Z(\widetilde{X}, \sO_{\widetilde{X}}(A)) \ar[d]^{F_{\widetilde{X}}^e}  \ar[r] & 0\\
0  \ar[r] &  H^{d-1}_{\widetilde{D} \cap Z}(\widetilde{D}, \sO_{\widetilde{D}}(qA)) \ar[r] & H^d_Z(\widetilde{X}, \sO_{\widetilde{X}}(qA-\widetilde{D})) \ar[r] & H^{d}_Z(\widetilde{X}, \sO_{\widetilde{X}}(qA)) \ar[r] & 0\\
}$$
\end{small}

\begin{cln}
The induced $e$-times Frobenius map 
$$F_{\widetilde{X}}^e:H^d_Z(\widetilde{X}, \sO_{\widetilde{X}}(A-\widetilde{D})) \to H^d_Z(\widetilde{X}, \sO_{\widetilde{X}}(qA-\widetilde{D}))$$
is injective for all $q=p^e$. 
\end{cln}

\begin{clnproof}
Since we are working in characteristic $p \gg 0$, by Lemma \ref{char p>>0}, the induced $e$-times Frobenius map $F_{\widetilde{X}}^e:H^d_Z(\widetilde{X}, \sO_{\widetilde{X}}(A)) \to H^d_Z(\widetilde{X}, \sO_{\widetilde{X}}(qA))$
is injective for all $q=p^e$.  
Hence, by the above commutative diagram, it suffices to show that the induced $e$-times Frobenius map
$$F_{\widetilde{D}}^e:H^{d-1}_{\widetilde{D} \cap Z}(\widetilde{D},\sO_{\widetilde{D}}(A)) \to  H^{d-1}_{\widetilde{D} \cap Z}(\widetilde{D}, \sO_{\widetilde{D}}(qA))$$
is injective for all $q=p^e$. 

We denote $B=\{\Delta\} +\epsilon \ \Div_X(c)$. Then $B$ has no common components  with the support of $D$.  
Let $\nu:D^{\nu} \to D$ be the normalization of $D$ and $g:\widetilde{D} \to D^{\nu}$ be the induced morphism.  
Then there exists an effective $\Q$-divisor $B^{\nu}$ on $D^{\nu}$, called the different of $B$ on $D^{\nu}$ (see \cite[\S 3]{Sh} for details), such that $K_{D^{\nu}}+B^{\nu}$ is $\Q$-Cartier and $K_{D^{\nu}}+B^{\nu}=\nu^*(K_X+D+B|_D)$. 
Since $H-tF$ is an $f$-ample $\Q$-divisor on $\widetilde{X}$ whose support has no common components with that of $\widetilde{D}$, 
$-A|_{\widetilde{D}}=-g^*(K_{D^{\nu}}+B^{\nu})+(H-tF)|_{\widetilde{D}}$ is a $g$-ample $\Q$-divisor on $\widetilde{D}$ whose fractional part has simple normal crossing support. 
It therefore follows from Lemma \ref{char p>>0} again that  $F_{\widetilde{D}}^e:H^{d-1}_{\widetilde{D} \cap Z}(\widetilde{D},\sO_{\widetilde{D}}(A)) \to  H^{d-1}_{\widetilde{D} \cap Z}(\widetilde{D}, \sO_{\widetilde{D}}(qA))$ is injective for all $q=p^e$. 
\end{clnproof}

From now on, we stick to the powers $q=p^e$ of $p$ such that $q\epsilon \geq m:=\lceil t n_0 \rceil+r+2s$ and $(q-1)(K_X+\Delta)$ is Cartier. Note that there exist infinitely many such $q=p^e$. 
By the definition of $A$, Claim 2 implies that the map
$$c^{m}F_{\widetilde{X}}^e:H^d_Z(\widetilde{X}, \sO_{\widetilde{X}}(A'-\widetilde{D})) \to H^d_Z(\widetilde{X}, \sO_{\widetilde{X}}(qA'-\widetilde{D}))$$
is injective for all such $q=p^e$, where $A'=f^*(K_X+\Delta)+tF-H$. 
Let $$\delta_e:H^d_{\m}(\omega_R^{(q)}((q-1)\Delta)) \to H^{d}_Z(\widetilde{X}, \sO_{\widetilde{X}}(qA'-\widetilde{D}))$$ be the Matlis dual of a natural inclusion map 
$$H^0(\widetilde{X}, \sO_{\widetilde{X}}(K_{\widetilde{X}}-\lfloor qA' \rfloor+\widetilde{D})) \hookrightarrow R( (1-q)(K_X+\Delta)).$$
Also, let $\delta:H^d_{\m}(\omega_R) \to H^{d}_Z(\widetilde{X}, \sO_{\widetilde{X}}(A'-\widetilde{D}))$ be the Matlis dual of a natural inclusion $\adj((X,\Delta);\a^t) \hookrightarrow R$.
Then we have the following commutative diagram with exact rows:
$$\xymatrix{
0 \ar[r] & \Ker(\delta) \ar[d] \ar[r] & H^d_{\m}(\omega_R) \ar[d]^{c^{m}F_R^e} \ar[r]^{\delta \qquad} & H^{d}_Z(\widetilde{X}, \sO_{\widetilde{X}}(A'-\widetilde{D})) \ar[d]^{c^{m}F_{\widetilde{X}}^e}  \ar[r]  & 0\\
0  \ar[r] & \Ker(\delta_e) \ar[r] & H^d_{\m}(\omega_R^{(q)}((q-1)\Delta)) \ar[r]^{\delta_e} & H^{d}_Z(\widetilde{X}, \sO_{\widetilde{X}}(qA'-\widetilde{D})) \ar[r] & 0\\
}$$
Here, $\Ker(\delta)$ (resp. $\Ker(\delta_e)$) can be considered as the annihilator in $H^d_{\m}(\omega_R)$ (resp. $H^d_{\m}(\omega_R^{(q)}((q-1)\Delta)$)) of $\adj((X,\Delta);\a^t)$ (resp. $H^0(\widetilde{X}, \sO_{\widetilde{X}}(\lceil G-qtF+qH \rceil))$ where $G=K_{\widetilde{X}}-f^*(K_X+\Delta)+\widetilde{D}$)
with respect to the pairing $R \times H^d_{\m}(\omega_R) \to  H^d_{\m}(\omega_R)$ (resp. $R \times H^d_{\m}(\omega_R^{(q)}((q-1)\Delta)) \to  H^d_{\m}(\omega_R^{(q)}((q-1)\Delta))$). 

Let $\xi \in H^d_{\m}(\omega_R) \setminus \Ker (\delta)$. 
By the above observation, one has $c^m\xi^q \notin \Ker(\delta_e)$, that is,
$$c^{m}\xi^q \notin \mathrm{Ann}_{H^d_{\m}(\omega_R^{(q)}((q-1)\Delta))}H^0(\widetilde{X}, \sO_{\widetilde{X}}(\lceil G-qtF+qH \rceil))$$
for all $q=p^e$ such that $q \epsilon \geq m$ and $(q-1)(K_X+\Delta)$ is Cartier.  
The multiplication by $c$
$$H^d_{\m}(\omega_R^{(q)}(\lceil q\Delta \rceil-D)) \xrightarrow{\times c} H^d_{\m}(\omega_R^{(q)}((q-1)\Delta))$$ 
forces $c^{m-1}\xi^q$ not to belong to the annihilator of $H^0(\widetilde{X}, \sO_{\widetilde{X}}(\lceil G-qtF+qH \rceil))$ in $H^d_{\m}(\omega_R^{(q)}(\lceil q\Delta \rceil-D))$.  
Then, by the definition of $n_0$, there exists an integer $0 \leq n \leq n_0$ such that 
$$c^{m-1}\xi^q \notin \mathrm{Ann}_{H^d_{\m}(\omega_R^{(q)}(\lceil q\Delta \rceil-D))}H^0(\widetilde{X}, \sO_{\widetilde{X}}((q-n)(H-tF))).$$
Since $H^0(\widetilde{X}, \sO_{\widetilde{X}}((q-n)(H-tF))) \subseteq \overline{\a^{\lceil t(q-n_0) \rceil}}$, we have $c^{m-1}\xi^q\overline{\a^{\lceil t(q-n_0) \rceil}} \ne 0$ in $H^d_{\m}(\omega_R^{(q)}(\lceil q\Delta \rceil-D))$. 
On the other hand, recall that $c \in \a$ and $m=\lceil t n_0 \rceil+r+2s$. 
Since $c^s \in \widetilde{\tau}^{\mathrm{div}}(R,D) \subseteq \widetilde{\tau}(R)$,  by the tight closure version of Brian\c con-Skoda's theorem (see \cite{HH1}), 
$$c^{m-s-1}\overline{\a^{\lceil t(q-n_0) \rceil}}=c^{\lceil tn_0 \rceil+r+s-1}\overline{\a^{\lceil t(q-n_0) \rceil}} \subseteq \widetilde{\tau}(R)\overline{\a^{\lceil tq \rceil+r-1}} \subseteq \a^{\lceil tq \rceil}.$$  
Thus, we have $c^{s} \a^{\lceil tq \rceil} \xi^q \neq 0$ in $H^d_{\m}(\omega_R^{(q)}(\lceil q \Delta \rceil-D))$ for all $q=p^e$ such that $q \epsilon \geq m$ and $(q-1)(K_X+\Delta)$ is Cartier.  
This implies that $\xi$ does not belong to $0_{H^d_{\m}(\omega_R)}^{\Div*(\Delta,\a^t)}$ by Proposition \ref{test}, because $c^{s}$ is a divisorial $D$-test element.  
So, we conclude that $0_{H^d_{\m}(\omega_R)}^{\Div*(\Delta, \a^t)} \subseteq \Ker(\delta)$, or dually $\widetilde{\tau}^{\mathrm{div}}((R,\Delta);\a^t) \supseteq \adj((X,\Delta);\a^t)$.
\end{proof}

\begin{cor}[\textup{\cite[Conjecture 5.1.1]{HW}}]\label{plt}
Let $(R,\m)$ be a normal local ring essentially of finite type over a field of characteristic zero, and let $\Delta$ be an effective $\R$-divisor on $X:=\Spec R$ such that $K_X+\Delta$ is $\R$-Cartier. 
Then $(X,\Delta)$ is plt if and only if $(R,\Delta)$ is of divisorially F-regular type. 
\end{cor}


\begin{thebibliography}{99}
 \bibitem{AM} Aberbach, I., MacCrimmon, B.: Some results on test elements. Proc. Edinb. Math. Soc., II. Ser. \textbf{42}, 541--549 (1999)

 \bibitem{A1} Aoyama, Y.: On the depth and the projective dimension of the canonical module. Japan. J. Math. (N.S.) \textbf{6}, 61--66  (1980)
 
 \bibitem{A2} Aoyama, Y.: Some basic results on canonical modules, J. Math. Kyoto Univ. \textbf{23}, 85--94 (1983)
 
 \bibitem{DH} Debarre, O.,  Hacon, C.\! D.: Singularities of divisors of low degree on abelian varieties. Manuscripta Math. \textbf{122} (2007), no. 2, 217--288.
 
 \bibitem{DEL} Demailly, J.-P., Ein, L., Lazarsfeld, R.: A subadditivity property of multiplier ideals. Mich. Math. J. \textbf{48}, 137--156  (2000) 
 
 \bibitem{EL} Ein, L., Lazarsfeld, R.: Singularities of theta divisors and the birational geometry of irregular varieties. J. Am. Math. Soc. \textbf{10}, 243--258  (1997) 
 
 \bibitem{Ha1} Hara, N.: A characterization of rational singularities in terms of injectivity of Frobenius map. Am. J. Math. \textbf{120}, 981--996 (1998)
 
 \bibitem{Ha2} Hara, N.: Geometric interpretation of tight closure and test ideals. Trans. Am. Math. Soc. \textbf{353}, 1885--1906 (2001)

 \bibitem{HT} Hara, N., Takagi, S.: On a generalization of test ideals. Nagoya Math. J. \textbf{175}, 59--74  (2004) 

 \bibitem{HW} Hara, N., Watanabe, K.-i.: F-regular and F-pure rings vs. log terminal and log canonical singularities. J. Algebr. Geom. \textbf{11}, 363--392  (2002) 

 \bibitem{HY} Hara, N., Yoshida, K.: A generalization of tight closure and multiplier ideals. Trans. Am. Math. Soc. \textbf{355}, 3143--3174  (2003)

\bibitem{Hi} Hironaka, H.: Resolution of singularities of an algebraic variety over a field of characteristic zero. I, II. Ann. Math. (2) \textbf{79}, 109--203 (1964); ibid. (2) \textbf{79}, 205--326  (1964)

 \bibitem{HH1} Hochster, M., Huneke, C.: Tight closure, invariant theory and the Brian\c con-Skoda theorem. J. Am. Math. Soc. \textbf{3}, 31--116  (1990)

 \bibitem{HH2} Hochster, M., Huneke, C.: Tight closure and strong F-regularity. Colloque en l'honneur de Pierre Samuel (Orsay, 1987). M\'em. Soc. Math. Fr., Nouv. S\'er. \textbf{38}, 119--133 (1989)
 
 \bibitem{Hu} Huneke, C.: Tight closure, parameter ideals, and geometry. Six lectures on commutative algebra (Bellaterra, 1996). Progr. Math. \textbf{166}, 187--239. Birkh\"auser, Basel (1998)
 
 \bibitem{Ka}  Kawakita, M.: Inversion of adjunction on log canonicity. Invent. Math. \textbf{167} (2007), no. 1, 129--133.
 
 \bibitem{KMM} Kawamata, Y., Matsuda, K., Matsuki, K.: Introduction to the minimal model problem. Algebraic geometry (Sendai, 1985).  Adv. Stud. Pure Math. \textbf{10}, 283--360. North-Holland, Amsterdam (1987) 
 
 \bibitem{Ko} Koll\'ar, J. et al.: Flips and abundance for algebraic threefolds. Ast\'erisque \textbf{211} (1992) 

 \bibitem{KM} Koll\'ar, J., Mori, S.: Birational geometry of algebraic varieties. Cambridge Tracts in Mathematics \textbf{134}, Cambridge University Press, Cambridge (1998) 
  
 \bibitem{La} Lazarsfeld, R.: Positivity in Algebraic Geometry II.  Ergebnisse der Mathematik und ihrer Grenzgebiete. 3. Folge, A Series of Modern Surveys in Mathematics, Vol. \textbf{49}, Springer-Verlag, Berlin (2004)
 
 \bibitem{MS} Mehta, V.B., Srinivas, V.: A characterization of rational singularities. Asian. J. Math. \textbf{1}, 249--278  (1997)

 \bibitem{Sh} Shokurov, V.V.: 3-fold log flips. Izv. Ross. Akad. Nauk, Ser. Mat. \textbf{56}, 105--203  (1992) 
 
 \bibitem{Sm1} Smith, K.E.: F-rational rings have rational singularities. Am. J. Math. \textbf{119}, 159--180   (1997)
 
 \bibitem{Sm} Smith, K.E.: The multiplier ideal is a universal test ideal. Commun. Algebra \textbf{28}, 5915--5929  (2000) 
 
 \bibitem{Ta} Takagi, S.: An interpretation of multiplier ideals via tight closure. J. Algebr. Geom. \textbf{13} , 393--415 (2004) 
 
 \bibitem{Ta3} Takagi, S.: Formulas for multiplier ideals on singular varieties. Am. J. Math. \textbf{128} (2006), no.6, 1345-1362.
\end{thebibliography}
\end{document}